\documentclass[10pt]{amsart}
\usepackage{latexsym, amsmath, amsfonts, amssymb, mathrsfs, fancyhdr, tikz, tikz-cd}
\usepackage{verbatim}
\usepackage{multicol}
\usetikzlibrary{matrix,arrows,decorations.pathmorphing,calc}

\def\ds{\displaystyle}

\def\n{\noindent}
\def\R{\mathbb{R}}

\def\v{\varphi}
\def\e{\varepsilon}
\def\s{\sigma}

\def\d{\delta}

\def\l{\lambda}

\def\V{\mathcal{V}}

\def\a{\alpha}
\def\b{\beta}

\def\g{\gamma}

\def\la{\langle}
\def\ra{\rangle}
\def\ae{a_\e}
\def\be{b_\e}
\def\ce{c_\e}
\def\de{d_\e}
\def\ee{e_\e}
\def\fe{f_\e}
\def\re{r_\e}

\def\ch2{\mathbb{C} \mathbb{H}^2}
\def\h2{\mathbb{H}^2}
\def\ddt{\frac{\partial}{\partial \theta}}
\def\ddr{\frac{\partial}{\partial r}}
\def\th{\theta}
\def\Rl{R_{\l}}
\def\oshr/2{\cosh \left( \frac{r}{2} \right)}
\def\inhr/2{\sinh \left( \frac{r}{2} \right)}
\def\osh2r/2{\cosh^2 \left( \frac{r}{2} \right)}
\def\inh2r/2{\sinh^2 \left( \frac{r}{2} \right)}
\def\-1/4{- \frac{1}{4}}
\def\H{\mathcal{H}}

\newtheorem{theorem}{Theorem}[section]
\newtheorem{lemma}[theorem]{Lemma}
\newtheorem{cor}[theorem]{Corollary}

\newtheorem{proposition}[theorem]{Proposition}

\theoremstyle{definition}

\theoremstyle{remark}
\newtheorem{remark}[theorem]{Remark}

\numberwithin{equation}{section}

\begin{document}

\title{Real hyperbolic hyperplane complements in the complex hyperbolic plane}

\author{Barry Minemyer}
\address{Department of Mathematics, The Ohio State University, Columbus, Ohio 43210}
\email{minemyer.1@osu.edu}



\date{\today.}



\begin{abstract}
This paper studies Riemannian manifolds of the form $M \setminus S$, where $M^4$ is a complete four dimensional Riemannian manifold with finite volume whose metric is modeled on the complex hyperbolic plane $\ch2$, and $S$ is a compact totally geodesic codimension two submanifold whose induced Riemannian metric is modeled on the real hyperbolic plane $\h2$.
In this paper we write the metric on $\ch2$ in polar coordinates about $S$, compute formulas for the components of the curvature tensor in terms of arbitrary warping functions (Theorem \ref{thm:curvature tensor}), and prove that there exist warping functions that yield a complete finite volume Riemannian metric on $M \setminus S$ whose sectional curvature is bounded above by a negative constant (Theorem \ref{thm:metric of negative curvature}(1)).  
The cases of $M \setminus S$ modeled on $\mathbb{H}^n \setminus \mathbb{H}^{n-2}$ and $\mathbb{C} \mathbb{H}^n \setminus \mathbb{C} \mathbb{H}^{n-1}$ were studied by Belegradek in \cite{Belegradek real} and \cite{Belegradek complex}, respectively.
One may consider this work as ``part 3" to this sequence of papers.
\end{abstract}

\maketitle



\section{Introduction}\label{Section:Introduction}
Let $M$ be a complete (connected) locally symmetric Riemannian manifold with finite volume and negative sectional curvature, and let $S$ be a (possibly disconnected) compact totally geodesic codimension two submanifold of $M$.  
It is known that the pair $(M,S)$ is modeled on $(\mathbb{H}^n, \mathbb{H}^{n-2})$, $(\mathbb{C} \mathbb{H}^n, \mathbb{C} \mathbb{H}^{n-1})$, or the ``exceptional case" $(\ch2, \h2)$.  
In \cite{Belegradek real} and \cite{Belegradek complex} Belegradek provides an in depth study of $M \setminus S$, the manifold obtained from $M$ by ``drilling out" $S$, when the model for the pair $(M,S)$ is one of the first two situations.
Here we consider the exceptional case, when $(M,S)$ is modeled on $(\ch2,\h2)$.  

The main result proved in this paper is the following.

\vskip 15pt

\begin{theorem}\label{thm:metric of negative curvature}
If $M$ is a complete finite volume complex hyperbolic $2$-manifold and $S$ is a compact totally real totally geodesic $2$-dimensional submanifold, then 
	\begin{enumerate}
	\item  $M \setminus S$ admits a complete finite volume metric with sectional curvature $\leq -1$.  
	\item  $M \setminus S$ admits a complete finite volume $A$-regular metric with sectional curvature $< 0$.
	\end{enumerate}
\end{theorem}

\vskip 10pt

The manifold $M \setminus S$ is diffeomorphic to a compact manifold $N$ obtained by cutting out a tubular neighborhood of $S$ in $M$ and removing all cusps in the case that $M$ is not compact.  
There are two possible types of boundary components of $N$.  
The first are compact infranil manifolds, which arise as cross-sections of the cusps of $M$ (if any), and the second type is circle bundles over the components of $S$.  


All three statements in the following Corollary \ref{cor:known corollary} can be deduced from other known results, as will be discussed in Remark \ref{rmk:corollaries}.
But combining Theorem \ref{thm:metric of negative curvature} with the methods deployed in \cite{Belegradek real} and \cite{Belegradek complex} provides independent proofs of these statements.

\vskip 15pt

\begin{cor}\label{cor:known corollary}
Suppose that $M$ is a complete finite volume complex hyperbolic $2$-manifold and $S$ is a compact totally real totally geodesic $2$-dimensional submanifold.  Then 
	\begin{enumerate}
	\item  the group $\pi_1(M \setminus S)$ is non-elementary (strongly) relatively hyperbolic, where the peripheral subgroups are the fundamental groups of the ends of $M \setminus S$.
	\item  $\pi_1(M \setminus S)$ satisfies the Farrell-Jones isomorphism conjecture. 
	\item  $\pi_1(M \setminus S)$ satisfies the Rapid Decay Property and the Baum-Connes conjecture.
	\end{enumerate}
\end{cor}

\vskip 10pt

\begin{remark}\label{rmk:corollaries}
The proof that Corollary \ref{cor:known corollary}(1) follows from Theorem \ref{thm:metric of negative curvature}(1) is identical to its analogues in \cite{Belegradek complex} and \cite{Belegradek real}.  
In particular, see Section 12 of \cite{Belegradek complex} in conjunction with Theorem 4.2 in \cite{Belegradek real}.
Many other properties of $\pi_1(M \setminus S)$ are then known to follow from Corollary \ref{cor:known corollary}(1) (see most conclusions of Theorem 1.1 in \cite{Belegradek real} and Theorem 1.4 in \cite{Belegradek complex}).  
But the fact that $\pi_1(M \setminus S)$ is relatively hyperbolic relative to the fundamental groups of its ends in our situation of $(M,S)$ being modeled on $(\ch2,\h2)$ is now a special case of Corollary 1.2 by Belegradek and Hruska in \cite{BH}.  

Results by Roushon (\cite{Roushon a} and \cite{Roushon b}) together with a deep result of Bartels-Farrell-L{\"u}ck (\cite{BFL}) imply that the fundamental groups of circle bundles over closed hyperbolic surfaces satisfy the Farrell-Jones isomorphism conjecture (\cite{FJ93}, abbreviated FJIC).  
A recent preprint by Bartels \cite{Bartels} proves that, if a countable group $G$ is relatively hyperbolic relative to a collection of subgroups that all satisfy the FJIC, then $G$ satisfies the FJIC.  
So combining Corollary \ref{cor:known corollary}(1) with Bartels' preprint proves Corollary \ref{cor:known corollary}(2).  


We can deduce the Rapid Decay Property (RDP) as follows.  
Drutu and Sapir \cite{DS} proved that if a finitely generated group is relatively hyperbolic relative to subgroups $P_1, \hdots, P_n$, then $G$ satisfies the RDP if and only if all of the subgroups $P_1, \hdots, P_n$ do. 
So it remains to check the RDP for the fundamental groups of the ends of $M \setminus S$.
But the ends of $M \setminus S$ are either infranil manifolds or circle bundles over components of $S$.  
In the first case, the RDP was established by Jolissaint in \cite{Jol90}, and in the second case it is established by Garncarek in \cite{Gar15} (also, one could combine \cite{Jol90} with \cite{Nos92}).
Finally, Lafforgue \cite{Lafforgue} proved that the fundamental group of any complete Riemannian manifold equipped with a non-positively curved $A$-regular metric which satisfies the RDP must also satisfy the Baum-Connes conjecture.  

\end{remark}

Theorem \ref{thm:metric of negative curvature}(1) is proved in Sections \ref{Section:metric in cylindrical coordinates} through \ref{Section:proof of negative curvature}, an outline of which is as follows.
Let $\h2$ denote a totally real totally geodesic $2$-plane in the complex hyperbolic plane $\ch2$.
We first express the metric in $\ch2$ in polar coordinates about $\h2$ (Section \ref{Section:metric in cylindrical coordinates}, or see below for the formula).  
We then, allowing for different coefficient functions in this metric, derive formulas for the components of the curvature tensor (Sections \ref{Section:basis and brackets} to \ref{Section:curvature tensor}, Theorem \ref{thm:curvature tensor}).
Our approach for Theorem \ref{thm:curvature tensor} is direct.  
We first fix a ``nice" non-holonomic frame of $\ch2 \setminus \h2$, compute the Lie brackets of these vector fields, compute the Levi-Civita connection with respect to this frame, and finally we compute the components of the curvature tensor.  
By far, the most difficult part of this is computing the Lie brackets (Theorem \ref{thm:bracket coefficients}).  
In Sections \ref{Section:h_th=h_r} and \ref{Section:sectional curvature} we prove a few general Lemma's about the sectional curvature functional.  
Finally, we construct specific ``warping" functions and prove that the sectional curvature of this warped metric is bounded above by a negative constant (Sections \ref{Section:construction of metric} and \ref{Section:proof of negative curvature}, which occupy the majority of this manuscript).  

Various types of warped metric computations have been used to prove some very important results in Riemannian geometry.  
Recently, they have been used by Ontaneda to prove the existence of ``smooth Riemannian hyperbolization" (\cite{Ontaneda hyperbolization}).  
A simple warped product was used by Gromov and Thurston in \cite{GT} in a key way, while more difficult warped product computations have been utilized by Farrell and Jones, or Farrell and Ontaneda, in a variety of papers (see \cite{Ontaneda} and the references therein).  
Therefore, the formulas in Theorem \ref{thm:curvature tensor} for the components of the curvature tensor of our warped metric may be the most useful portion of this paper.  
So let us describe them.

Let $r$ denote the distance from a point to the totally real totally geodesic $2$-plane $\h2$, and let $\ddr$ be the unit length vector field on $\ch2 \setminus \h2$ pointing radially from $\h2$.  
Then there exists an orthogonal, non-holonomic set of vector fields $(S,T)$ on $\ch2 \setminus \h2$ such that the complex hyperbolic metric $g$ is
	\begin{equation*}
	g = \osh2r/2 dS^2 + \cosh^2(r) dT^2 + \inh2r/2 d \theta^2 + dr^2.
	\end{equation*}
In the above formula, $dS$ and $dT$ denote the covector fields dual to $S$ and $T$, respectively, and $d \theta$ denotes the standard metric on the unit circle.  
See Sections \ref{Section:metric in cylindrical coordinates} through \ref{Section:curvature tensor} for more details, especially about $S$ and $T$.
Let
	\begin{equation*}
	\l = h_\th^2(r) dS^2 + h_r^2(r) dT^2 + v^2(r) d \theta^2 + dr^2
	\end{equation*}
where $h_\th, h_r$, and $v$ are positive functions of $r$.  
Let $Y_1 = \frac{1}{v} \ddt$, $Y_2 = \frac{1}{h_\th} S$, $Y_3 = \frac{1}{h_r} T$, and $Y_4 = \ddr$ be a $\l$-orthonormal frame of $\ch2 \setminus \h2$.  
Then the formulas for the components of the curvature tensor of $\l$ with respect to $h_\th$, $h_r$, and $v$ are

\vskip 10pt

	\begin{align*}
	&\la \Rl(Y_1,Y_2)Y_1,Y_2 \ra_{\l} = - \frac{v'h_\th'}{v h_\th} - \frac{1}{4} \left( \frac{-v^2}{4h_\th^2h_r^2} - \frac{h_\th^2}{4v^2h_r^2} + \frac{3h_r^2}{4v^2h_\th^2} - \frac{1}{2v^2} + \frac{1}{2h_\th^2} - \frac{1}{2h_r^2} \right)  \\
	&\la \Rl(Y_1,Y_3)Y_1,Y_3 \ra_{\l} = - \frac{v'h_r'}{vh_r} - \frac{1}{4} \left( \frac{- v^2}{4h_\th^2h_r^2} + \frac{3h_\th^2}{4v^2h_r^2} - \frac{h_r^2}{4v^2h_\th^2} - \frac{1}{2v^2} - \frac{1}{2h_\th^2} + \frac{1}{2h_r^2}  \right) \\
	&\la \Rl(Y_2,Y_3)Y_2,Y_3 \ra_{\l} = - \frac{h_\th' h_r'}{h_\th h_r} - \frac{1}{4} \left( \frac{3v^2}{4h_\th^2 h_r^2} - \frac{h_\th^2}{4v^2h_r^2} - \frac{h_r^2}{4v^2h_\th^2} + \frac{1}{2v^2} + \frac{1}{2h_\th^2} + \frac{1}{2h_r^2} \right)  \\
	&\la \Rl(Y_1,Y_2)Y_3,Y_4 \ra_{\l} = \frac{-1}{4h_r} \left[ \left( \frac{h_\th}{v} \right)' - \left( \frac{v}{h_\th} \right)' - \left( \frac{h_r^2}{vh_\th} \right)' \right]  \\
	&\la \Rl(Y_1,Y_3)Y_2,Y_4 \ra_{\l} = \frac{-1}{4h_\th} \left[ - \left( \frac{h_r}{v} \right)' + \left( \frac{v}{h_r} \right)' + \left( \frac{h_\th^2}{vh_r} \right)' \right]  
	\end{align*}
	\begin{align*}
	&\la \Rl(Y_1,Y_4)Y_2,Y_3 \ra_{\l} = \frac{-1}{4v} \left[ \left( \frac{h_\th}{h_r} \right)' + \left( \frac{h_r}{h_\th} \right)' + \left( \frac{v^2}{h_\th h_r} \right)' \right] \\
	&\la \Rl(Y_1,Y_4)Y_1,Y_4 \ra_{\l} = - \frac{v''}{v} \hskip 80pt \la \Rl(Y_2,Y_4)Y_2,Y_4 \ra_{\l} = - \frac{h_\th''}{h_\th} \\
	&\la \Rl(Y_3,Y_4)Y_3,Y_4 \ra_{\l} = - \frac{h_r''}{h_r}
	\end{align*}
	
\vskip 10pt
	
\n and where all other components of the curvature tensor are identically zero.  
Of course, the formulas above intimately depend on the vector fields $S$ and $T$ (which are parallel to $J \ddt$ and $J \ddr$, respectively).  
These vector fields are constructed in Section \ref{Section:metric in cylindrical coordinates}.

We end this paper by proving Theorem \ref{thm:metric of negative curvature}(2) in Section \ref{Section:A-regular metric}.  
This construction and proof are nearly identical to Section 11 of \cite{Belegradek complex}, and we only include some details here due to the slight differences in both our metric and curvature formulas (see Section \ref{Section:h_th=h_r}).
In fact, there are several places in this paper where we refer to results in  \cite{Belegradek complex}, and others where we follow \cite{Belegradek complex} very closely.  
So the reader interested in understanding every detail of these results should also have a copy of \cite{Belegradek complex} at hand.
In fact, one may consider this paper as ``part 3" in the sequence of papers \cite{Belegradek real} and \cite{Belegradek complex}, or maybe more accurately just the sequel to \cite{Belegradek complex}.  
Many of the situations that arise in this work are considerably different than those in \cite{Belegradek complex}, but whenever possible we have tried to use results from Belegradek's paper in order to simplify calculations here.  
Lastly, we have in most cases tried to use the same notation as \cite{Belegradek complex} in order to make it easier to simultaneously read the two papers.

\begin{remark}\label{rmk:curvature notation}
One major notational difference between this paper and the papers \cite{Belegradek complex} and \cite{Belegradek real} is the following.
Let $g$ be a Riemannian metric with Levi-Civita connection $\nabla$, and let $W, X, Y,$ and $Z$ be vector fields.  
In this paper we use the notation
	\begin{equation*}
	R(X,Y)Z = \nabla_Y \nabla_X Z - \nabla_X \nabla_Y Z + \nabla_{[X,Y]} Z
	\end{equation*}
for the curvature tensor $R$ of $g$.  
The negative of this formula is used in \cite{Belegradek complex} and \cite{Belegradek real}.
So, in particular, the $(4,0)$-curvature tensor $\la R(X,Y)Z,W \ra_g$ in this paper is equivalent to $\la R(X,Y)W,Z \ra_g$ in \cite{Belegradek complex} and \cite{Belegradek real}.  
\end{remark}

\vskip 20pt

\section{The metric in complex hyperbolic 2-space in cylindrical coordinates about a totally real totally geodesic 2-plane}\label{Section:metric in cylindrical coordinates}

The purpose of this section is to describe the metric for $\ch2$ in cylindrical coordinates about a totally real totally geodesic $2$-plane denoted by $\h2$.  
The case of $\mathbb{C}\mathbb{H}^n$ about $\mathbb{C}\mathbb{H}^{n-1}$ was worked out in \cite{Belegradek complex}, and the metric on distance spheres in $\mathbb{C}\mathbb{H}^n$  is discussed by Farrell and Jones in \cite{Farrell Jones}.  
The current case is worked out by Phillips in an undergraduate REU \cite{Phillips}, but in the following sections we will rely on the vector fields $S$ and $T$ defined below in a very crucial way.  
So what follows is a detailed explanation of the construction in \cite{Phillips}.
The terminology and notation in this section will closely follow \cite{Belegradek complex} and \cite{Belegradek real}.  

Let $g$ denote the complex hyperbolic metric on $\ch2$ normalized to have constant holomorphic sectional curvature $-1$.  
Since $\h2$ is a complete totally geodesic submanifold of the negatively curved manifold $\ch2$, there exists an orthogonal projection map $\pi : \ch2 \to \h2$.  
This map $\pi$ is a fiber bundle whose fibers are totally real totally geodesic $2$-planes, and therefore have constant sectional curvature $\-1/4$.  

For $r > 0$ let $F(r)$ denote the $r$-neighborhood of $\h2$.  
Then $F(r)$ is a real hypersurface in $\ch2$, and consequently we can decompose $g$ as
	\begin{equation*}
	g = dr^2 + g_r
	\end{equation*}
where $g_r$ is the induced Riemannian metric on $F(r)$.  
Let $\pi_r : F(r) \to \h2$ denote the restriction of $\pi$ to $F(r)$.  
Note that $\pi_r$ is a circle bundle whose fiber over any point $q \in \ch2$ is the circle of radius $r$ in the totally real totally geodesic $2$-plane $\pi^{-1}(q)$.  
The tangent bundle splits as an orthogonal sum $\mathcal{V}(r) \oplus \mathcal{H}(r)$ where $\mathcal{V}(r)$ is tangent to the circle $\pi_r^{-1}(q)$ and $\mathcal{H}(r)$ is the orthogonal complement to $\mathcal{V}(r)$.

For $r, v > 0$ there exists a diffeomorphism $\phi_{vr}: F(v) \to F(r)$ induced by the geodesic flow along the totally real totally geodesic $2$-planes orthogonal to $\h2$.   
Fix $p \in F(r)$ arbitrary, let $q = \pi(p) \in \h2$, and let $\g$ be the unit speed geodesic such that $\g(0) = q$ and $\g(r) = p$. 
In what follows, all computations are considered in the tangent space $T_p F(r)$. 

Note that $\V(r)$ is tangent to both $F(r)$ and the totally real totally geodesic $2$-plane $\pi^{-1}(q)$.  
Then since $\pi^{-1}(q)$ is preserved by the geodesic flow, we have that $d \phi_{vr}$ takes $\V(v)$ to $\V(r)$.  
Since $exp_p^{-1} \left( \pi^{-1}(q) \right)$ is a real 2-plane with curvature $\-1/4$, the metric $g$ restricted to $\pi^{-1}(q)$ can be written as $dr^2 + \inh2r/2 d \th^2$ where $d \th^2$ is the standard metric on the unit circle.
 Note that the vector field $\ddt$ is invariant under $d \phi_{vr}$.

Let $J$ denote the complex structure on $\ch2$.  
It is well known that $J_p$ preserves complex lines in $T_p \ch2$ and maps real $2$-planes into (and, via dimension reasons in our setting, onto) their orthogonal complement.
Since $\left( \ddt, \ddr \right)$ spans a real $2$-plane in $T_p \ch2$, their orthogonal complement $\H_p(r)$ is spanned by $\left( J \ddt, J \ddr \right)$.  
In the following two Subsections we define vector fields $S$ and $T$ which are just scaled copies of $J \ddt$ and $J \ddr$, respectively.

\vskip 20pt
\subsection{Definition of the vector field S}

First note that $\left( \ddr, J \ddt \right)$ spans a real $2$-plane in $T_p \ch2$ (since its $J$-image is its orthogonal complement).  
So $P = \text{exp} \left( \text{span} \left( \ddr, J \ddt \right) \right)$ is a totally real totally geodesic $2$-plane in $\ch2$ which intersects $\h2$ orthogonally.
Since this intersection is orthogonal, $P$ is preserved by the geodesic flow $\phi$.  
Therefore, $\text{span} \left( J \ddt \right)$ is preserved by $d \phi$.  

The set $P \cap \h2$ is a (real) geodesic.  
Let $\a(s)$ denote this geodesic parameterized with respect to arc length so that $\a(0) = q$.  
Then define $S_p = (d \pi)^{-1}_p \a'(0)$.  
There exists a positive real-valued function $a(r,s)$ so that the metric $g$ restricted to $P$ is of the form $dr^2 + a^2(r,s) ds^2$.  
But note that $a(r,s)$ is independent of $s$ because of the isometric $\R$-action on $P$ by translations along $\a$.  
Then since the curvature of a real $2$-plane is $\-1/4$, we have that $a(r) = \oshr/2$.

\vskip 20pt
\subsection{Definition of the vector field T}

This is analogous to the definition of $S$.
Note that $\left( \ddr, J \ddr \right)$ spans a complex line in $T_p \ch2$ (since it is preserved by its $J$-image).
So $Q = \text{exp} \left( \text{span} \left( \ddr, J \ddr \right) \right)$ is a complex geodesic in $\ch2$ which intersects $\h2$ orthogonally.
Since this intersection is orthogonal, $Q$ is preserved by the geodesic flow $\phi$.  
Therefore, $\text{span} \left( J \ddr \right)$ is preserved by $d \phi$.  

The set $Q \cap \h2$ is a (real) geodesic.  
Let $\b(t)$ denote this geodesic parameterized with respect to arc length so that $\b(0) = q$.
Then define $T_p = (d \pi)^{-1} \b'(0)$.  
There exists a positive real-valued function $b(r,t)$ so that the metric $g$ restricted to $Q$ is of the form $dr^2 + b^2(r,t) dt^2$.
But note that $b(r,t)$ is independent of $t$ because of the isometric $\R$-action on $Q$ by translations along $\b$. 
Then since the curvature of a complex geodesic is $-1$, we have that $b(r) = \cosh r$.

\vskip 20pt
\subsection{Conclusion}


\begin{theorem}
The complex hyperbolic manifold $\ch2 \setminus \h2$ can be written as $(0, \infty) \times F$  where $F \cong \R^2 \times \mathbb{S}^1$ equipped with the metric
	\begin{equation}\label{eqn: complex hyperbolic metric}
	g = \inh2r/2 d \th^2 + \osh2r/2 dS^2 + \cosh^2 (r) \, dT^2 + dr^2 .
	\end{equation}
\end{theorem}


In equation \eqref{eqn: complex hyperbolic metric}, $dS$ and $dT$ denote the covector fields dual to the vector fields $S$ and $T$, respectively.
Lastly, notice that $dS^2 + dT^2$ is the hyperbolic metric with constant sectional curvature $\-1/4$.

\vskip 20pt

\section{The metric $\l$ and Lie brackets (Part I)}\label{Section:basis and brackets}

Fix an open interval $I$ and let $v$, $h_{\th}$, and $h_r$ be smooth positive functions on I.  
Let $\l$ denote the Riemannian metric
	 \begin{equation*}
	 \l = v^2 d \th^2 + h_{\th}^2 dS^2 + h_r^2 dT^2 + dr^2
	 \end{equation*}
on $I \times F$ (and where $F = F(r)$ for some generic $r > 0$).
Of course, we recover the metric on $\ch2 \setminus \h2$ when $v = \inhr/2$, $h_\th = \oshr/2$, and $h_r = \cosh r$.  
The purpose of Sections \ref{Section:basis and brackets} through \ref{Section:curvature tensor} is to compute the components of the curvature tensor $R_\l$ in terms of $v$, $h_\l$, and $h_r$.    

For these curvature computations we will use the non-holonomic basis
	\begin{equation}\label{eqn:basis}
	X_1 = \ddt, \qquad X_2 = S, \qquad X_3 = T, \qquad X_4 = \ddr.
	\end{equation}
Let us note the following observations about this basis:
	\begin{enumerate}
	\item $\la X_1, X_1 \ra_\l = v^2$, $\la X_2, X_2 \ra_\l = h_\th^2$, $\la X_3, X_3 \ra_\l = h_r^2$, and $\la X_4, X_4 \ra_\l = 1$.   
	\item $[X_i, X_4] = 0$ since each $X_i$ is invariant under the flow of $\ddr$.  
	\item $\la [X_i, X_j],X_4 \ra_\l = 0$ since $[X_i,X_j]$ is tangent to level surfaces of $r$.  
	\end{enumerate}
\vskip 10pt
\n By (3) above, there exist constants $(a_i), (b_j), $ and $(c_k)$, $1 \leq i, j, k \leq 3$, such that
	\begin{align}
	[X_1,X_2] &= a_1 X_1 + a_2 X_2 + a_3 X_3  \label{eqn:x1x2} \\
	[X_1,X_3] &= b_1 X_1 + b_2 X_2 + b_3 X_3 \label{eqn:x1x3} \\
	[X_2,X_3] &= c_1 X_1 + c_2 X_2 + c_3 X_3  \label{eqn:x2x3}.
	\end{align}
The purpose of Section \ref{Section:Lie brackets} is to compute these constants.

The orthonormal basis corresponding to \eqref{eqn:basis} is
	\begin{equation}\label{eqn:orthonormal basis}
	Y_1 = \frac{1}{v}X_1, \qquad Y_2 = \frac{1}{h_\th}X_2, \qquad Y_3 = \frac{1}{h_r}X_3, \qquad Y_4 = X_4.
	\end{equation}
Direct calculations show that the Lie brackets for this basis have the following properties:
	\begin{enumerate}
	\item $[Y_1,Y_2] = \frac{1}{v h_\th}[X_1,X_2]$, $[Y_1,Y_3] = \frac{1}{v h_r}[X_1,X_3]$, and $[Y_2, Y_3] = \frac{1}{h_\th h_r}[X_2,X_3]$.  
	\item $[Y_1,Y_4] = \frac{v'}{v}Y_1$, $[Y_2,Y_4] = \frac{h_\th'}{h_\th}Y_2$, and $[Y_3,Y_4]=\frac{h_r'}{h_r}Y_3$.
	\end{enumerate}

\vskip 20pt

\section{Components of the curvature tensor in $\ch2$}\label{Section:curvature in ch2}

The components of the (4,0) curvature tensor of the complex hyperbolic metric $g$ can be expressed in terms of $g$ and the complex structure $J$.  
The following formula can be found in \cite{KN} or in Section 5 of \cite{Belegradek complex} (recall Remark \ref{rmk:curvature notation} from the Introduction).  
In this formula $X, Y, Z$, and $W$ are arbitrary vector fields.
	\begin{align*}
	4 & \la R_g(X,Y)Z,W \ra_g = \la X,W \ra_g \la Y,Z \ra_g - \la X,Z \ra_g \la Y,W \ra_g  \\
	&+  \la X,JW \ra_g \la Y,JZ \ra_g - \la X,JZ \ra_g \la Y,JW \ra_g + 2 \la X,JY \ra_g \la W,JZ \ra_g.
	\end{align*}
Recall that the complex structure $J$ on $\ch2$ preserves the complex hyperbolic metric $g$.  
That is, for any vector fields $X$ and $Y$, $\la X,Y \ra_g = \la JX, JY \ra_g$.  
We therefore have that
	\begin{equation*}
	JY_1 = Y_2, \qquad JY_2 = -Y_1, \qquad JY_3 = -Y_4, \qquad JY_4 = Y_3.
	\end{equation*}

Up to symmetries of the curvature tensor, the following are the only non-zero components of the curvature tensor of the complex hyperbolic metric $g$ with respect to the orthonormal basis $(Y_1, Y_2, Y_3, Y_4)$.

	\begin{multicols}{2}
	\n \begin{equation} \la R_g(Y_1,Y_2)Y_1,Y_2 \ra_g = - 1 \label{eqn:curvature 1}  \end{equation}
	\begin{equation} \la R_g(Y_1,Y_3)Y_1,Y_3 \ra_g = \-1/4 \label{eqn:curvature 2}  \end{equation}
	\begin{equation} \la R_g(Y_2,Y_3)Y_2,Y_3 \ra_g = - \frac{1}{4} \label{eqn:curvature 3}  \end{equation}
	\begin{equation} \la R_g(Y_1,Y_2)Y_3,Y_4 \ra_g = \frac{1}{2} \label{eqn:curvature 4}  \end{equation}
	\begin{equation} \la R_g(Y_1,Y_3)Y_2,Y_4 \ra_g = \frac{1}{4} \label{eqn:curvature 5}  \end{equation}
	\begin{equation} \la R_g(Y_1,Y_4)Y_2,Y_3 \ra_g = - \frac{1}{4} \label{eqn:curvature 6}  \end{equation}
	\begin{equation} \la R_g(Y_1,Y_4)Y_1,Y_4 \ra_g = - \frac{1}{4} \label{eqn:curvature 7}  \end{equation}
	\begin{equation} \la R_g(Y_2,Y_4)Y_2,Y_4 \ra_g = \-1/4 \label{eqn:curvature 8}  \end{equation}
	\begin{equation} \la R_g(Y_3,Y_4)Y_3,Y_4 \ra_g = - 1 \label{eqn:curvature 9}  \end{equation}
	\end{multicols}
Equations \eqref{eqn:curvature 1} through \eqref{eqn:curvature 9} will be used in the following Section to determine the values of the coefficients in equations \eqref{eqn:x1x2} through \eqref{eqn:x2x3}.

\vskip 20pt

\section{Lie brackets (Part II)}\label{Section:Lie brackets}

The purpose of this Section is to compute the coefficients in equations \eqref{eqn:x1x2}, \eqref{eqn:x1x3}, and \eqref{eqn:x2x3}.
The major tool is a formula worked out by Belegradek in \cite{Belegradek real} and stated in Appendix B of \cite{Belegradek complex}.  
The set up for this formula is as follows.  
Suppose that $\l = dr^2 + \l_r$ is a warped product metric on $I \times F$, and let $\{ X_i \}$ be a $\l_r$-orthogonal basis of vector fields on a neighborhood $U \subset F$.  
Let $h_i(r) = \sqrt{\la X_i, X_i \ra_{\l_r}}$ so that the collection $Y_i = \frac{1}{h_i} X_i$ forms an orthonormal basis on $U$ for any $r>0$.  
Then (recall Remark \ref{rmk:curvature notation} for the difference between this and \cite{Belegradek complex})
	\begin{align}
	2&\la R \left( \ddr, Y_i \right) Y_j, Y_k \ra_\l =   \label{eqn:Belegradek} \\
	\la [Y_i,Y_k],Y_j \ra_\l \left( \ln \frac{h_j}{h_k} \right)' &+ \la [Y_j,Y_i],Y_k \ra_\l \left( \ln \frac{h_k}{h_j} \right)' + \la [Y_j,Y_k],Y_i \ra_\l \left( \ln \frac{h_i^2}{h_j h_k} \right)' . \nonumber
	\end{align}

As a first (easy) step, we use equation \eqref{eqn:Belegradek} to prove the following.

\vskip 10pt

\begin{lemma}\label{lemma:zero coefficients}
In equations \eqref{eqn:x1x2}, \eqref{eqn:x1x3}, and \eqref{eqn:x2x3}
	\begin{equation}\label{eqn:zero coefficients}
	a_1 = a_2 = b_1 = b_3 = c_2 = c_3 = 0.
	\end{equation}
\end{lemma}

\vskip 4pt

\begin{proof}
We will prove that $a_1 = 0$, and then indicate how to analogously show that each of the other coefficients in equation \eqref{eqn:zero coefficients} are zero.

With respect to the complex hyperbolic metric $g$, we know that $\la R_g(Y_4,Y_1)Y_1,Y_2 \ra_g = 0$.  
But then by equation \eqref{eqn:Belegradek} we have that
	\begin{align*}
	0 &= \la [Y_1,Y_2],Y_1 \ra_g \left( \ln \frac{h_v}{h_\th} \right)' + 0 + \la [Y_1,Y_2],Y_1 \ra_g \left( \ln \frac{h_v}{h_\th} \right)'  \\
	&= \frac{2}{\inh2r/2 \oshr/2} \la [X_1,X_2],X_1 \ra_g \left( \ln \left( \frac{\inhr/2}{\oshr/2} \right) \right)'  \\
	&= \frac{2 a_1}{\oshr/2} \left( \ln \tanh \left( \frac{r}{2} \right) \right)'
	\end{align*}
and thus $a_1 = 0$.  

Then to show:
	\begin{align*}
	&a_2 = 0 \text{ use } \la R_g(Y_4,Y_2)Y_1,Y_2 \ra_g = 0  \\
	&b_1 = 0 \text{ use } \la R_g(Y_4,Y_1)Y_1,Y_3 \ra_g = 0  \\
	&b_3 = 0 \text{ use } \la R_g(Y_4,Y_3)Y_1,Y_3 \ra_g = 0  \\
	&c_2 = 0 \text{ use } \la R_g(Y_4,Y_2)Y_2,Y_3 \ra_g = 0  \\
	&c_3 = 0 \text{ use } \la R_g(Y_4,Y_3)Y_2,Y_3 \ra_g = 0 .
	\end{align*}
\end{proof}

Due to Lemma \ref{lemma:zero coefficients}, we can rewrite equations \eqref{eqn:x1x2}, \eqref{eqn:x1x3}, and \eqref{eqn:x2x3} as
	\begin{align}
	&[X_1,X_2] = \a_3 X_3 \label{eqn:x3}  \\
	&[X_1,X_3] = \a_2 X_2 \label{eqn:x2}  \\
	&[X_2,X_3] = \a_1 X_1 \label{eqn:x1}
	\end{align}
where $\a_3 = a_3$, $\a_2 = b_2$, and $\a_1 = c_1$.  
The reason that the $\a_i$ are not zero is because of equations \eqref{eqn:curvature 4}, \eqref{eqn:curvature 5}, and \eqref{eqn:curvature 6}.  
If these ``mixed components" of the curvature tensor for $g$ were zero, then these coefficients would also be zero.  
In particular, this would be the case if we were dealing with $\mathbb{H}^4$ instead of $\ch2$.  

But we can combine equations \eqref{eqn:curvature 4} through \eqref{eqn:curvature 6} with \eqref{eqn:Belegradek}, and use the {\it Nijenhuis tensor}, to prove the following Theorem.

\vskip 20pt

\begin{theorem}\label{thm:bracket coefficients}
The coefficients in equations \eqref{eqn:x3}, \eqref{eqn:x2}, and \eqref{eqn:x1} are
	\begin{equation}
	\a_1 = \frac{1}{2}, \qquad \a_2 = \frac{1}{2}, \qquad \a_3 = - \frac{1}{2}.
	\end{equation}
\end{theorem}

\vskip 20pt

\begin{proof}
First, note that
	\begin{equation}
	[Y_1,Y_2] = \frac{\a_3 h_r}{v h_\th} \, Y_3, \qquad [Y_1,Y_3] = \frac{\a_2 h_\th}{v h_r} \, Y_2, \qquad [Y_2,Y_3] = \frac{\a_1 v}{h_\th h_r} \, Y_1
	\end{equation}
where, since we are dealing with the complex hyperbolic metric $g$, $v = \inhr/2$, $h_\th = \oshr/2$, and $h_r = \cosh(r)$.  

We first combine equations \eqref{eqn:curvature 5} and \eqref{eqn:curvature 6} with equation \eqref{eqn:Belegradek} to obtain
	\begin{align}
	- \frac{1}{2} = \la [Y_2,&Y_3],Y_1 \ra_g \left( \ln \frac{v}{h_r} \right)' + \la [Y_1,Y_2],Y_3 \ra_g \left( \ln \frac{h_r}{v} \right)' + \la [Y_1,Y_3],Y_2 \ra_g \left( \ln \frac{h_\th^2}{v h_r} \right)'  \nonumber \\
	&= \frac{a_1 v}{h_\th h_r} \left( \ln \frac{v}{h_r} \right)' - \frac{a_3 h_r}{v h_\th} \left( \ln \frac{v}{h_r} \right)' + \frac{a_2 h_\th}{v h_r} \left( \ln \frac{h_\th^2}{v h_r} \right)'  \nonumber \\
	&= \left( \frac{\a_1 \inhr/2}{\oshr/2 \cosh(r)} - \frac{\a_3 \cosh(r)}{\inhr/2 \oshr/2} \right) \left( \frac{1}{2} \coth \left( \frac{r}{2} \right) - \tanh(r) \right)  \label{eqn:equation 1} \\
	&+ \frac{\a_2 \oshr/2}{\inhr/2 \cosh(r)} \left( \tanh \left( \frac{r}{2} \right)  - \frac{1}{2} \coth \left( \frac{r}{2} \right) - \tanh(r)  \right) \nonumber
	\end{align}
and
	\begin{align}
	- \frac{1}{2} = \la [Y_1,&Y_2],Y_3 \ra_g \left( \ln \frac{h_r}{h_\th} \right)' + \la [Y_3,Y_1],Y_2 \ra_g \left( \ln \frac{h_\th}{h_r} \right)' + \la [Y_3,Y_2],Y_1 \ra_g \left( \ln \frac{v^2}{h_\th h_r} \right)'  \nonumber \\
	&= \frac{a_3 h_r}{v h_\th} \left( \ln \frac{h_r}{h_\th} \right)' + \frac{a_2 h_\th}{v h_r} \left( \ln \frac{h_r}{h_\th} \right)' - \frac{a_1 v}{h_\th h_r} \left( \ln \frac{v^2}{h_\th h_r} \right)'  \nonumber \\
	&= \left( \frac{\a_3 \cosh(r)}{\inhr/2 \oshr/2} + \frac{\a_2 \oshr/2}{\inhr/2 \cosh(r)} \right) \left( \tanh(r) - \frac{1}{2} \tanh \left( \frac{r}{2} \right) \right)  \label{eqn:equation 2} \\
	&- \frac{\a_1 \inhr/2}{\oshr/2 \cosh(r)} \left( \coth \left( \frac{r}{2} \right) - \frac{1}{2} \tanh \left( \frac{r}{2} \right) - \tanh(r) \right).  \nonumber
	\end{align}

This yields two equations with three unknown variables.  
It is a tedious excercise to check that the values in Theorem \ref{thm:bracket coefficients} satisfy these two equations.
But, of course, there is no reason to believe (yet) that this solution is unique.
In an attempt to obtain a third independent equation we could combine equation \eqref{eqn:curvature 4} with equation \eqref{eqn:Belegradek}, but one easily checks that this leads to a dependent system of equations.

There are two routes to obtaining a third independent equation involving the $\a_i's$.  
One way is to compute the components of $R_\l$ as equations involving the $\a_i's$ and $v$, $h_\th$, and $h_r$.  
We could then substitute in the values for $g$ to obtain several other independent equations.  
But the easier method is to use the {\it Nijenhuis Tensor}.  
Since the almost complex structure on $\ch2$ is integrable, we have that 
	\begin{equation}\label{eqn:Nijenhuis}
	0 = [X,Y] + J[JX,Y] + J[X,JY] - [JX,JY]
	\end{equation}
for all vector fields $X$ and $Y$.  
Evaluating \eqref{eqn:Nijenhuis} at $X = Y_1$ and $Y=Y_3$ yields
	\begin{align}
	0 &= [Y_1,Y_3] + J[Y_2,Y_3] + J[Y_1,-Y_4] - [Y_2,-Y_4]  \nonumber \\
	&= \frac{\a_2 h_\th}{v h_r} Y_2 + \frac{\a_1 v}{h_\th h_r} JY_1 - \frac{v'}{v}JY_1 + \frac{h_\th'}{h_\th}Y_2 \nonumber \\
	&= \left( \frac{\a_2 \oshr/2}{\inhr/2 \cosh(r)} + \frac{\a_1 \cosh(r)}{\oshr/2 \inhr/2} - \frac{1}{2} \coth \left( \frac{r}{2} \right) + \frac{1}{2} \tanh \left( \frac{r}{2} \right) \right) Y_2  \nonumber  \\
	&\Longrightarrow \qquad \frac{1}{\cosh(r)} \left( \coth \left( \frac{r}{2} \right) \a_2 + \tanh \left( \frac{r}{2} \right) \a_1 \right) = \frac{1}{\sinh(r)}.  \label{eqn:equation 3}
	\end{align}
Now one easily checks that equations \eqref{eqn:equation 1}, \eqref{eqn:equation 2}, and \eqref{eqn:equation 3} are independent, guaranteeing a unique solution.
It is also easy this time to check that the values in Theorem \ref{thm:bracket coefficients} satisfy equation \eqref{eqn:equation 3}, completing the proof of Theorem \ref{thm:bracket coefficients}.

\end{proof}

The following lemma just restates what has been proven so far:

\vskip 15pt

\begin{lemma}\label{lem:Lie brackets}
The six independent Lie brackets for the vector fields described in equation \eqref{eqn:orthonormal basis} are
	\begin{align*}
	[Y_1,&Y_2] = \frac{-h_r}{2vh_\th} \,  Y_3, \qquad [Y_1,Y_3] = \frac{h_\th}{2v h_r} \, Y_2, \qquad [Y_2,Y_3] = \frac{v}{2 h_\th h_r} \, Y_1  \\
	&[Y_1,Y_4] = \frac{v'}{v} \, Y_1, \qquad [Y_2,Y_4] = \frac{h_\th'}{h_\th} \, Y_2, \qquad [Y_3,Y_4] =  \frac{h_r'}{h_r} \, Y_3.
	\end{align*}
\end{lemma}

\vskip 20pt

\section{The Levi-Civita connection of $\l$}\label{Section:connection}

In this Subsection we will compute the Levi-Civita connection $\nabla$ associated to the metric $\l$ with respect to the frame $( Y_1, Y_2, Y_3, Y_4 )$.  
To perform this calculation we will use the well-known ``Koszul formula" (which can be found on pg. 55 of \cite{do Carmo})
	\begin{align}
	\la \nabla_Y X, Z \ra_\l = & \, \frac{1}{2} ( X \la Y,Z \ra_\l + Y \la Z,X \ra_\l - Z \la X,Y \ra_\l  \label{eqn:Kozul formula}  \\
	&- \la [X,Z],Y \ra_\l - \la [Y,Z], X \ra_\l - \la [X,Y],Z \ra_\l ) . \nonumber
	\end{align}
Since we are considering an orthonormal frame we know that $\la Y_i , Y_j \ra_\l = \d_{ij}$, where $\d_{ij}$ denotes Kronecker's delta.
Therefore the first three terms on the right hand side of formula \eqref{eqn:Kozul formula} are all zero in our given frame.  
Thus, formula \eqref{eqn:Kozul formula} reduces to
	\begin{equation}
	\la \nabla_Y X, Z \ra_\l = - \frac{1}{2} \left( \la [X,Z],Y \ra_\l + \la [Y,Z], X \ra_\l + \la [X,Y],Z \ra_\l \right).  \label{eqn:connection}
	\end{equation}
	
It is now a simple calculation using formula \eqref{eqn:connection} and the results of Lemma \ref{lem:Lie brackets} to prove the following Theorem.

\vskip 15pt

\begin{theorem}\label{thm:connection}
The Levi-Civita connection $\nabla$ is given by the $16$ equations
	\begin{align}
	&0 = \nabla_{Y_4}Y_1 = \nabla_{Y_4}Y_2 = \nabla_{Y_4}Y_3 = \nabla_{Y_4}Y_4 \label{eqn61} \\
	\nabla_{Y_1}Y_1 = - &\frac{v'}{v} Y_4, \qquad \nabla_{Y_2}Y_2 = - \frac{h_\th'}{h_\th} Y_4, \qquad \nabla_{Y_3}Y_3 = - \frac{h_r'}{h_r} Y_4  \label{eqn62}\\
	\nabla_{Y_1}Y_4 = &\frac{v'}{v} Y_1, \qquad \nabla_{Y_2}Y_4 = \frac{h_\th'}{h_\th} Y_2, \qquad \nabla_{Y_3}Y_4 = \frac{h_r'}{h_r} Y_3  \label{eqn63}  \\
	&\nabla_{Y_1}Y_2 = - \frac{1}{2} \left( \frac{v}{2 h_\th h_r} + \frac{h_\th}{2 v h_r} + \frac{h_r}{2 v h_\th} \right) Y_3  \label{eqn64} \\
	&\nabla_{Y_1}Y_3 = - \frac{1}{2} \left( \frac{-v}{2 h_\th h_r} - \frac{h_\th}{2 v h_r} - \frac{h_r}{2 v h_\th} \right) Y_2  \label{eqn65} \\
	&\nabla_{Y_2}Y_1 = - \frac{1}{2} \left( \frac{v}{2 h_\th h_r} + \frac{h_\th}{2 v h_r} - \frac{h_r}{2 v h_\th} \right) Y_3  \label{eqn66} \\
	&\nabla_{Y_2}Y_3 = - \frac{1}{2} \left( \frac{-v}{2 h_\th h_r} - \frac{h_\th}{2 v h_r} + \frac{h_r}{2 v h_\th} \right) Y_1  \label{eqn67}  \\
	&\nabla_{Y_3}Y_1 = - \frac{1}{2} \left( \frac{-v}{2 h_\th h_r} + \frac{h_\th}{2 v h_r} - \frac{h_r}{2 v h_\th} \right) Y_2  \label{eqn68}  \\
	&\nabla_{Y_3}Y_2 = - \frac{1}{2} \left( \frac{v}{2 h_\th h_r} - \frac{h_\th}{2 v h_r} + \frac{h_r}{2 v h_\th} \right) Y_1  \label{eqn69} 
	\end{align}
\end{theorem}

\vskip 10pt

\begin{proof}
We only prove equation \eqref{eqn64} and the first equalities in equations \eqref{eqn61}, \eqref{eqn62}, and \eqref{eqn63}.  
All of the other equations are proven analogously.  
As stated above, each equation is obtained by simply combining Lemma \ref{lem:Lie brackets} with formula \eqref{eqn:connection}.
To prove the first equality in \eqref{eqn61} we compute
	\begin{equation*}
	\la \nabla_{Y_4}Y_1, Y_i \ra_\l = - \frac{1}{2} \left( \la [Y_1,Y_i],Y_4 \ra_\l + \la [Y_4,Y_i], Y_1 \ra_\l + \la [Y_1,Y_4],Y_i \ra_\l \right).
	\end{equation*}
Using Lemma \ref{lem:Lie brackets} we see that each term above is zero when $i = 2, 3,$ or $4$.  
When $i=1$, the first term is zero and the last two terms cancel due to the antisymmetry of the Lie bracket.  
Thus, $\nabla_{Y_4}Y_1 = 0$.  

To prove the first equality in \eqref{eqn62} we plug into formula \eqref{eqn:connection} to obtain
	\begin{align*}
	\la \nabla_{Y_1}Y_1, Y_i \ra_\l &= - \frac{1}{2} \left( \la [Y_1,Y_i],Y_1 \ra_\l + \la [Y_1,Y_i], Y_1 \ra_\l + \la [Y_1,Y_1],Y_i \ra_\l \right)  \\
	&= - \la [Y_1,Y_i],Y_1 \ra_\l .
	\end{align*}
This is nonzero only when $i=4$, and substituting the value of $[Y_1,Y_4]$ given in Lemma \ref{lem:Lie brackets} yields the desired result.

For the first equality in \eqref{eqn63} we have
	\begin{equation*}
	\la \nabla_{Y_1}Y_4, Y_i \ra_\l = - \frac{1}{2} \left( \la [Y_4,Y_i],Y_1 \ra_\l + \la [Y_1,Y_i], Y_4 \ra_\l + \la [Y_4,Y_1],Y_i \ra_\l \right).
	\end{equation*}
Each term is zero when $i = 2, 3$ or $4$.  
When $i=1$ we have 
	\begin{equation*}
	\la \nabla_{Y_1}Y_4, Y_1 \ra_\l = - \la [Y_4,Y_1],Y_1 \ra_\l = - \left( - \frac{v'}{v} \la Y_1,Y_1 \ra_\l \right) = \frac{v'}{v}.
	\end{equation*}
	
Lastly, to verify equation \eqref{eqn64}, we have from formula \eqref{eqn:connection}
	\begin{equation*}
	\la \nabla_{Y_1}Y_2, Y_i \ra_\l = - \frac{1}{2} \left( \la [Y_2,Y_i],Y_1 \ra_\l + \la [Y_1,Y_i], Y_2 \ra_\l + \la [Y_2,Y_1],Y_i \ra_\l \right).
	\end{equation*}
We see from Lemma \ref{lem:Lie brackets} that each term above is zero when $i = 1, 2,$ or $4$.  
When $i=3$ we have that
	\begin{align*}
	\la \nabla_{Y_1}Y_2, Y_3 \ra_\l &= - \frac{1}{2} \left( \la [Y_2,Y_3],Y_1 \ra_\l + \la [Y_1,Y_3], Y_2 \ra_\l - \la [Y_1,Y_2],Y_3 \ra_\l \right)  \\
	&= - \frac{1}{2} \left( \frac{v}{2 h_\th h_r} + \frac{h_\th}{2 v h_r} - \frac{- h_r}{2 v h_\th} \right).
	\end{align*}
\end{proof}

\vskip 20pt

\section{Components of the curvature tensor of $\l$}\label{Section:curvature tensor}

Recall from Remark \ref{rmk:curvature notation} that in this paper we are using the definition
	\begin{equation}\label{eqn:def of R}
	R(X,Y)Z = \nabla_Y \nabla_X Z - \nabla_X \nabla_Y Z + \nabla_{[X,Y]} Z.
	\end{equation}
for the curvature tensor.
Let $\Rl$ denote the curvature tensor of the metric $\l$.  
Then, up to symmetries of the curvature tensor, the only non-zero components of $\Rl$ are given by the following Theorem.

\vskip 20pt

\begin{theorem}\label{thm:curvature tensor}
In terms of the basis given in equation \eqref{eqn:orthonormal basis}, the only independent nonzero components of the $(4,0)$ curvature tensor $\Rl$ are the following:
	\begin{align*}
	&\la \Rl(Y_1,Y_2)Y_1,Y_2 \ra_{\l} = - \frac{v'h_\th'}{v h_\th} - \frac{1}{4} \left( \frac{-v^2}{4h_\th^2h_r^2} - \frac{h_\th^2}{4v^2h_r^2} + \frac{3h_r^2}{4v^2h_\th^2} - \frac{1}{2v^2} + \frac{1}{2h_\th^2} - \frac{1}{2h_r^2} \right)  \\
	&\la \Rl(Y_1,Y_3)Y_1,Y_3 \ra_{\l} = - \frac{v'h_r'}{vh_r} - \frac{1}{4} \left( \frac{- v^2}{4h_\th^2h_r^2} + \frac{3h_\th^2}{4v^2h_r^2} - \frac{h_r^2}{4v^2h_\th^2} - \frac{1}{2v^2} - \frac{1}{2h_\th^2} + \frac{1}{2h_r^2}  \right) \\
	&\la \Rl(Y_2,Y_3)Y_2,Y_3 \ra_{\l} = - \frac{h_\th' h_r'}{h_\th h_r} - \frac{1}{4} \left( \frac{3v^2}{4h_\th^2 h_r^2} - \frac{h_\th^2}{4v^2h_r^2} - \frac{h_r^2}{4v^2h_\th^2} + \frac{1}{2v^2} + \frac{1}{2h_\th^2} + \frac{1}{2h_r^2} \right)  \\
	&\la \Rl(Y_1,Y_2)Y_3,Y_4 \ra_{\l} = \frac{-1}{4h_r} \left[ \left( \frac{h_\th}{v} \right)' - \left( \frac{v}{h_\th} \right)' - \left( \frac{h_r^2}{vh_\th} \right)' \right]  \\
	&\la \Rl(Y_1,Y_3)Y_2,Y_4 \ra_{\l} = \frac{-1}{4h_\th} \left[ - \left( \frac{h_r}{v} \right)' + \left( \frac{v}{h_r} \right)' + \left( \frac{h_\th^2}{vh_r} \right)' \right] \\
	&\la \Rl(Y_1,Y_4)Y_2,Y_3 \ra_{\l} = \frac{-1}{4v} \left[ \left( \frac{h_\th}{h_r} \right)' + \left( \frac{h_r}{h_\th} \right)' + \left( \frac{v^2}{h_\th h_r} \right)' \right] \\
	&\la \Rl(Y_1,Y_4)Y_1,Y_4 \ra_{\l} = - \frac{v''}{v} \hskip 80pt \la \Rl(Y_2,Y_4)Y_2,Y_4 \ra_{\l} = - \frac{h_\th''}{h_\th} \\
	&\la \Rl(Y_3,Y_4)Y_3,Y_4 \ra_{\l} = - \frac{h_r''}{h_r}.
	\end{align*}
\end{theorem}
	
\vskip 10pt
	
It is an exercise in hyperbolic trigonometric identities to check that, when $v = \inhr/2$, $h_\th = \oshr/2$, and $h_r = \cosh(r)$, the above formulas reduce to the constants in equations \eqref{eqn:curvature 1} through \eqref{eqn:curvature 9}.  
Also, note that the curvature tensor on a four dimensional Riemannian manifold has $21$ components which are independent with respect to the symmetries of the curvature tensor.  
So Theorem \ref{thm:curvature tensor} also states that the remaining $12$ components of $\Rl$ are identically zero.

\begin{proof}
There is nothing enlightening about computing the $21$ independent components of $\Rl$.  
We simply plug the results of Lemma \ref{lem:Lie brackets} and Theorem \ref{thm:connection} into equation \eqref{eqn:def of R} and punch out the computations using the properties of $\nabla$ (and remembering that $Y_4 = \ddr$ and all of $v$, $h_\th$, and $h_r$ are functions of $r$).  
As an illustration, let us prove the first equality in Theorem \ref{thm:curvature tensor}.  
In what follows, $a = v/(2h_\th h_r)$, $b = h_\th/(2v h_r)$, and $c = h_r/(2 v h_\th)$.  
	\begin{align*}
	\la \Rl(Y_1,Y_2)Y_1,Y_2 \ra_{\l} &= \la \nabla_{Y_2} \nabla_{Y_1} Y_1 - \nabla_{Y_1} \nabla_{Y_2} Y_1 + \nabla_{[Y_1,Y_2]} Y_1, Y_2 \ra_\l  \\
	&= \la -\frac{v'}{v} \nabla_{Y_2} Y_4 + \frac{1}{2} (a+b-c) \nabla_{Y_1}Y_3 - c \nabla_{Y_3}Y_1, Y_2 \ra_\l  \\
	&= \la - \frac{v' h_\th'}{v h_\th} Y_2 - \frac{1}{4} (a+b-c)(-a-b-c) Y_2 + \frac{1}{2} c(-a+b-c) Y_2, Y_2 \ra_\l  \\
	&= - \frac{v' h_\th'}{v h_\th} - \frac{1}{4} (-a^2 - b^2 + 3c^2 - 2ab + 2ac - 2bc)  \\
	&= - \frac{v'h_\th'}{v h_\th} - \frac{1}{4} \left( \frac{-v^2}{4h_\th^2h_r^2} - \frac{h_\th^2}{4v^2h_r^2} + \frac{3h_r^2}{4v^2h_\th^2} - \frac{1}{2v^2} + \frac{1}{2h_\th^2} - \frac{1}{2h_r^2} \right).
	\end{align*}
\end{proof}

\vskip 20pt

\section{The case when $h_\th = h_r$}\label{Section:h_th=h_r}

In Sections \ref{Section:construction of metric} and \ref{Section:proof of negative curvature} we will construction functions for $v$, $h_\th$, and $h_r$ for which $\l$ will be complete, have finite volume, and have negative curvature bounded away from zero.  
The purpose of Sections \ref{Section:h_th=h_r} and \ref{Section:sectional curvature} are to derive formulas and results which will help to prove that the metric constructed in Section \ref{Section:construction of metric} has sectional curvature bounded above by a negative constant.
The formulas in Theorem \ref{thm:curvature tensor} are rather long and complicated, but in the very special case when $h_\th = h_r := h$ these formulas reduce very nicely as stated in the following Corollary.
In what follows, $K(Y_i, Y_j) := \la R(Y_i,Y_j)Y_i,Y_j \ra_\l$ denotes the sectional curvature (with respect to $\l$) of the 2-plane spanned by $(Y_i,Y_j)$.  

\vskip 15pt

\begin{cor}[Corollary to Theorem \ref{thm:curvature tensor}]\label{cor:curvature tensor}
When $h_\th = h_r := h$, the formulas in Theorem \ref{thm:curvature tensor} reduce to
	\begin{align*}
	&K(Y_1,Y_2) = K(Y_1,Y_3) = - \frac{v' h'}{vh} + \frac{v^2}{16h^4}  \\
	&K(Y_2,Y_3) = - \frac{1}{4h^2} - \frac{3v^2}{16h^4} - \left( \frac{h'}{h} \right)^2  \\
	&K(Y_1,Y_4) = - \frac{v''}{v} \hskip 80pt K(Y_2,Y_4) = K(Y_3,Y_4) = - \frac{h''}{h}  
	\end{align*}
	\begin{align*}
	&\la \Rl(Y_1,Y_4)Y_2,Y_3 \ra_{\l} = \frac{-1}{4v} \left( \frac{v^2}{h^2} \right)' = \frac{-v}{2h^2} \left( \ln \frac{v}{h} \right)'  \\
	&\la \Rl(Y_1,Y_2)Y_3,Y_4 \ra_{\l} = \frac{v}{4h^2} \left( \ln \frac{v}{h} \right)'  = - \la \Rl(Y_1,Y_3)Y_2,Y_4 \ra_{\l}
	\end{align*}
\end{cor}

\vskip 10pt

\begin{remark}
It is important to note that the complex hyperbolic metric $g$ on $\ch2 \setminus \h2$ does {\it not} satisfy the conditions of Corollary \ref{cor:curvature tensor}.  
Namely, $h_\th = \oshr/2 \neq \cosh(r) = h_r$.  
But we will use the above formulas in Sections \ref{Section:construction of metric} and \ref{Section:proof of negative curvature} in order to greatly reduce calculations.
\end{remark}

Comparing the equations in Corollary \ref{cor:curvature tensor} with equations (9.2) through (9.5) in \cite{Belegradek complex}, one sees that they are nearly identical.  
The only difference is that the above equations contain specific (and inconsistent) values for the constant $c_{23}$ in \cite{Belegradek complex}.  
More specifically, the constant $c_{23}$ takes on the values of $0$ and $1/2$ in the equation for $K(Y_2,Y_3)$ (from left to right), and takes on the value of $1/2$ in the formula for $\la \Rl(Y_1,Y_2)Y_3,Y_4 \ra_{\l}$.

In order to simplify calculations even further, let us prove the following Lemma before computing a formula for the sectional curvature of a generic $2$-plane when $h_\th = h_r$.

\vskip 15pt

\begin{lemma}\label{lemma:choosing basis}
Assume that $h_\th = h_r := h$, and suppose that $(U_2, U_3)$ is a orthonormal set of vectors (with respect to $\l$) whose span is the plane spanned by $(Y_2, Y_3)$.  
Furthermore, assume that $(U_2, U_3)$ and $(Y_2, Y_3)$ have the same orientation.  
Then the curvature formulas in Corollary \ref{cor:curvature tensor} remain unchanged if $Y_i$ is replaced with $U_i$ for $i = 2,3$.  
\end{lemma}

\vskip 10pt

\begin{proof}

By the assumptions, there exists constants $a_2, a_3, b_2, b_3$ such that
	\begin{align*}
	U_2 = a_2 Y_2 + &a_3 Y_3,	\qquad	U_3 = b_2 Y_2 + b_3 Y_3  \\
	a_2^2 &+ a_3^2 = 1 = b_2^2 + b_3^2   \\
	&a_2 b_2 + a_3 b_3 = 0 
	\end{align*}
	
Now notice that
	\begin{align*}
	(a_2 b_3 - a_3 b_2)^2 &= (a_2 b_3 - a_3 b_2)^2 + (a_2 b_2 + a_3 b_3)^2  \\
	&= (a_2^2 + a_3^2)(b_2^2 + b_3^2)  \\
	&= 1
	\end{align*}
and therefore $a_2b_3 - a_3 b_2 = \pm 1$.  
But since $(U_2, U_3)$ and $(Y_2, Y_3)$ have the same orientation, we have that $a_2b_3 - a_3 b_2 = 1$.  

Now we just use this formula in conjunction with the formulas in Corollary \ref{cor:curvature tensor}:
	\begin{align*}
	K(Y_1,U_2) = \la R(Y_1,U_2)Y_1,U_2 \ra &= a_2^2 \la R(Y_1,Y_2)Y_1,Y_2 \ra + a_3^2 \la R(Y_1,Y_3)Y_1,Y_3 \ra  \\
	&= (a_2^2 + a_3^2) \left( - \frac{v'h'}{vh} + \frac{v^2}{16h^4} \right) = - \frac{v'h'}{vh} + \frac{v^2}{16h^4}.
	\end{align*}
The proof that $K(Y_1,U_3) = K(Y_1,Y_3)$ is completely analogous, as is the proof that $K(U_2,Y_4) = K(U_3,Y_4) = K(Y_2,Y_4) = K(Y_3,Y_4)$.
Also, it is clear that $K(U_2,U_3) = K(Y_2,Y_3)$.  
So all that is left is to check the ``cross-terms" of the curvature tensor.  
Each of these three cases are nearly identical, so we only compute one here.
	\begin{align*}
	\la R(Y_1,U_2)U_3,Y_4 \ra &= \la R(Y_1,a_2Y_2 + a_3Y_3) b_2Y_2 + b_3Y_3, Y_4 \ra  \\
	&= a_2 b_3 \la R(Y_1,Y_2)Y_3,Y_4 \ra + a_3 b_2 \la R(Y_1,Y_3)Y_2,Y_4 \ra  \\
	&= (a_2 b_3 - a_3 b_2) \la R(Y_1,Y_2)Y_3,Y_4 \ra  \\
	&= \la R(Y_1,Y_2)Y_3,Y_4 \ra
	\end{align*}
and where the third equality above is due to the specific formulas in Corollary \ref{cor:curvature tensor}.
	
\end{proof}

Using the notation of Section \ref{Section:metric in cylindrical coordinates}, let $p \in (0,\infty) \times F$ and let $q = \pi(p) \in \h2$.  
Let $\sigma$ denote a generic $2$-plane tangent to $(0,\infty) \times F$ at $p$.  
In the generic case when $d\pi(\s) = \h2$, an identical argument to that given by Belegradek in (\cite{Belegradek complex}, Section 9, pg. 559) shows that there exists an orthonormal basis $(C,D)$ of $\s$ such that
	\begin{equation*}
	C = c_1 Y_1 + c_2 U_2 + c_3 U_3 + c_4 Y_4, \qquad D = d_1 Y_1 + d_2 U_2
	\end{equation*}
and $(U_2, U_3)$ are an orthonormal set as in Lemma \ref{lemma:choosing basis}.  

In the following calculation we use some new notation.
Since we proved in Lemma \ref{lemma:choosing basis} the curvature formulas are identical, in what follows we replace $(U_2,U_3)$ with $(Y_2,Y_3)$, respectively.
Also, we use the convention $R_{ijkl} := \la R(Y_i,Y_j)Y_k,Y_l \ra_\l$.  

We now compute:
	\begin{align*}
	K(\s) = K&(C,D) = d_1^2 K(C,Y_1) + d_2^2 K(C,Y_2) + 2 d_1 d_2 \la R(C,Y_1)C,Y_2 \ra  \\
	&K(C,Y_1) = c_2^2 K(Y_1,Y_2) + c_3^2 K(Y_1,Y_3) + c_4^2 K(Y_1,Y_4)  \\
	&K(C,Y_2) = c_1^2 K(Y_1,Y_2) + c_3^2 K(Y_2,Y_3) + c_4^2 K(Y_2,Y_4).
	\end{align*}
Also,
	\begin{align*}
	\la R(C,Y_1)C,Y_2 \ra &= -c_1 c_2 K(Y_1,Y_2) + c_3 c_4 R_{1324} + c_3 c_4 R_{1423}  \\
	&= -c_1 c_2 K(Y_1,Y_2) + \frac{3}{2} c_3 c_4 R_{1423}
	\end{align*}
since $R_{1324} = \frac{1}{2} R_{1423}$ whenever $h_\th = h_r$.  

Putting this all together gives that
	\begin{align}
	K(\s) &= (c_1d_2 - c_2d_1)^2 K(Y_1,Y_2) + d_1^2c_3^2K(Y_1,Y_3) + d_1^2c_4^2K(Y_1,Y_4) \label{eqn:curvature formula} \\
	&+ d_2^2c_3^2K(Y_2,Y_3) + d_2^2c_4^2K(Y_2,Y_4) + 3c_3c_4d_1d_2R_{1423}. \nonumber
	\end{align}
	
In an identical manner as Remark 9.6 in \cite{Belegradek complex}, if $R_{1423} = 0$ and $K(Y_i,Y_j)$ is less than or equal to a negative constant for each $i \neq j$ then $K(\s)$ is bounded above by the same negative constant.  
This is because the coefficients of the sectional curvatures of the coordinate planes in equation \eqref{eqn:curvature formula} add up to one.

Lastly, the case when $d \pi(\s) \neq \h2$ is identical to Remark 9.7 in \cite{Belegradek complex}.  
If dim($d \pi(\s)) = 0$, then $K(\s) = K(Y_1,Y_4)$.  
If dim($d \pi(\s)) = 1$, then we can choose our orthonormal basis $(C,D)$ in such a way that $K(\s)$ contains no ``mixed terms" (such as $R_{1423}$ in equation \eqref{eqn:curvature formula}).  
Therefore in the case when dim($d \pi(\s)) \neq 2$, if the sectional curvatures of the coordinate planes are bounded above by a negative constant, so is $K(\s)$.

\vskip 20pt

\section{Sectional curvature of $\l$}\label{Section:sectional curvature}
The purpose of this section is to set-up and prove Lemma \ref{lemma:general curvature formula}, Lemma \ref{lemma:continuity lemma}, and Corollary \ref{cor:final step} below, which will help us deal with the general case when $h_\th \neq h_r$.  
As above, let $p \in (0,\infty) \times F$ and let $\sigma$ denote a $2$-plane tangent to $(0,\infty) \times F$ at $p$.  
We can always find an orthonormal basis $(A,B)$ of $\s$ such that
	\begin{equation*}
	A = a_1Y_1 + a_2Y_2 + a_3Y_3 + a_4Y_4, \qquad B = b_2Y_2 + b_3Y_3 + b_4Y_4.
	\end{equation*}
Then, using notation from Section \ref{Section:h_th=h_r}, we have that
	\begin{align}
	&K(\s) = \la R(A,B)A,B \ra_\l  \label{eqn:general curvature formula}  \\
	= a_1^2 b_2^2 R_{1212} + a_1^2b_3^2 &R_{1313} + a_1^2b_4^2R_{1414} + (a_2b_3 - a_3b_2)^2 R_{2323} + (a_2b_4 - a_4b_2)^2 R_{2424}  \nonumber  \\
	+ (a_3b_4 - a_4&b_3)^2R_{3434} + 2a_1b_2(a_3b_4 - a_4b_3)R_{1234}   \nonumber  \\
	+ 2a_1b_3(a_2b_4& - a_4b_2)R_{1324} + 2a_1b_4(a_2b_3 - a_3b_2)R_{1423}.  \nonumber
	\end{align}
Equation \eqref{eqn:general curvature formula} is used to prove the following.

\vskip 15pt

\begin{lemma}\label{lemma:general curvature formula}
Suppose that $R_{1234}, R_{1324}, R_{1423} \neq 0$, and 
	\begin{align*}
	&(1a) \hskip 10pt R_{1212} \leq - | R_{1234} | \hskip 80pt (1b) \hskip 10pt R_{3434} \leq - | R_{1234} |  \\
	&(2a) \hskip 10pt R_{1313} \leq - | R_{1324} | \hskip 80pt (2b) \hskip 10pt R_{2424} \leq - | R_{1324} |  \\
	&(3a) \hskip 10pt R_{1414} \leq - | R_{1423} | \hskip 80pt (3b) \hskip 10pt R_{2323} \leq - | R_{1423} |  .
	\end{align*}
Then there exists $C < 0$ such that $K(\s) < C$ for any $2$-plane $\s$.
\end{lemma}

\vskip 10pt

\begin{remark}
Obviously, Lemma \ref{lemma:general curvature formula} also applies in the special case when $h_\th = h_r$.
\end{remark}

\begin{remark}
Also, it is ok if some of $R_{1234}, R_{1324},$ or $R_{1423} = 0$.  
But then for the conclusion of the Lemma to be true we need the inequalities using these components with value $0$ to be strict.
\end{remark}

\begin{proof}
With a little bit of arithmetic, one can rewrite equation \eqref{eqn:general curvature formula} for $K(\s)$ as:
	\begin{align*}
	-  (\pm a_1b_2 &+ a_3b_4 - a_4b_3)^2 |R_{1234}|  \\
	&+ a_1^2b_2^2(R_{1212} + |R_{1234}|) + (a_3b_4 - a_4b_3)^2 (R_{3434} + |R_{1234}|)  \\
	- (\pm a_1b_3 &+ a_2b_4 - a_4b_2)^2 |R_{1324}|  \\
	&+ a_1^2b_3^2(R_{1313} + |R_{1324}|) + (a_2b_4 - a_4b_2)^2 (R_{2424} + |R_{1324}|)  \\
	- (\pm a_1b_4 &+ a_2b_3 - a_3b_2)^2 |R_{1423}|  \\
	&+ a_1^2b_4^2(R_{1414} + |R_{1423}|) + (a_2b_3 - a_3b_2)^2 (R_{1414} + |R_{1423}|)
	\end{align*}
where the signs of the three ``$\pm$" terms depend on the signs of $R_{1234}$, $R_{1324}$, and $R_{1423}$, respectively.  
For example, if $R_{1234} < 0$ then the first term would be ``$- a_1 b_2 $", and similarly for $R_{1234} > 0$ and for the other two terms.  

Due to the inequalities in the assumptions of the Lemma, each term above is nonpositive.
A tedious but mostly trivial algebraic argument (which makes great use of the fact that $(A,B)$ is orthonormal) can be used to show that the coefficients of the first, fourth, and seventh terms cannot all simultaneously be $0$.  
Thus $K(\s) < 0$ for each fixed $\s$.  

Let $M$ denote the maximal value attained by $K()$ over all orthonormal pairs $(A,B)$.  
By compactness $M$ is attained via some $K(\s)$, and due to the last paragraph this maximum is strictly less than $0$.  
Thus, $K(\s) \leq M < 0$ for all $2$-planes $\s \in T_pM$.  
	
\end{proof}

The following ``Continuity Lemma" is a direct application of the equation in the proceeding proof.  
This Lemma will be used to deal with the endpoints of the regions in Sections \ref{Section:construction of metric} and \ref{Section:proof of negative curvature}, as well as to prove Corollary \ref{cor:final step}.

\vskip 15pt

\begin{lemma}[Continuity Lemma]\label{lemma:continuity lemma}
Let $\l$, $\l'$ be Riemannian metrics on a manifold $M$, and fix a point $p \in M$.  
Let $R$, $R'$ and $K$, $K'$ denote the curvature tensors and sectional curvatures of $\l$ and $\l'$, respectively.  
Then for all $\e > 0$ there exists $\d > 0$ such that if $|R_{ijkl} - R_{ijkl}' | < \d$ for all $i,j,k,l$ (at $p$), then $|K(\s) - K'(\s)| < \e$ for all $\s \in T_p M$.
\end{lemma}

\vskip 15pt

\begin{cor}\label{cor:final step}
Suppose that $R_{1234}, R_{1324}, R_{1423} \neq 0$, and that five of the six inequalities assumed in Lemma \ref{lemma:general curvature formula} hold.  
Then there exists $\d > 0$ such that, if the sixth inequality fails by at most $\d$, then $K()$ is still bounded above by a negative constant.
\end{cor}

\vskip 10pt

\begin{proof}
If the sixth inequality is instead a strict equality, then by Lemma \ref{lemma:general curvature formula} $K()$ is bounded above by a negative constant.  
Then just apply Lemma \ref{lemma:continuity lemma} to the two components of the curvature tensor that appear in that sixth inequality.
\end{proof}

\n We will make great use of Corollary \ref{cor:final step} in the final step in Section \ref{Section:proof of negative curvature}.  


\vskip 20pt

\section{Constructing a complete negatively curved metric with finite volume}\label{Section:construction of metric}

Our goal is to construct functions $v$, $h_\th$, and $h_r$ so that the metric $\l$ is complete, has finite volume, and so that the sectional curvature of $\l$ is bounded above by a negative constant.  
In this Section we simply construct these functions.  
In Section \ref{Section:proof of negative curvature} we prove the desired properties of the metric.  
It will be clear that the metric is complete and has finite volume (or see Remark 3.3 of \cite{Belegradek real}), so we will prove that the sectional curvature is bounded above by a negative constant.
The metric $\l$ will also agree with the complex hyperbolic metric $g$ when $r$ is at least half of the normal injectivity radius of the compact totally geodesic submanifold $S$.  
In this construction, the domain of $v$, $h_\th$, and $h_r$ will be from negative infinity to the normal injectivity radius of $S$, turning $S$ into a cusp of $M \setminus S$.  

There are seven stages to this process, as illustrated in Figure \ref{figure:warping stages}.  
The endpoints of the regions for the stages will all depend on a small positive constant $\e$ and will be defined below.  

In the first region, whose domain is $(- \infty, \ae)$ where $\ae < 0$, we set $v = \e e^r$ and $h_\th = h_r = e^{\frac{r}{2}}$.  
In region $2$, whose domain is $(\ae,\be)$ with $\be > 0$, we simultaneously ``warp" both $h_\th$ and $h_r$ from $e^{\frac{r}{2}}$ to $\oshr/2$.
Then in region $3$, defined on $(\be, \ce)$, we simply have that $v = \e e^r$ and $h_\th = h_r = \oshr/2$.  
Region $4$ occurs over the interval $(\ce,\de)$.  
Here we ``bend" $v$ from $\e e^r$ to $\inhr/2$, while keeping $h_\th = h_r = \oshr/2$.  
So on region $5$, occurring from $(\de, \ee)$, we have $v = \inhr/2$ and $h_\th = h_r = \oshr/2$.  

Notice that in these first five regions we have been very careful to keep $h_\th = h_r$ so that we can use the much simpler formulas in both Corollary \ref{cor:curvature tensor} and equation \eqref{eqn:curvature formula}.  
Since these formulas are very similar to those obtained in \cite{Belegradek complex}, much of this work can be copied over in order to compute the curvature bounds in these regions (and our warping functions will essentially agree with Belegradek's).  
Region $7$, defined from $\fe$ to the normal injectivity radius of $S$, is simple since here $\l$ will agree with the complex hyperbolic metric $g$ and thus the curvature will be bounded above by $- \frac{1}{4}$.  
But in region $6$, defined on $(\ee, \fe)$, we warp $h_r$ from $\oshr/2$ to $\cosh(r)$.  
So in this region we have $h_\th \neq h_r$, forcing us to use the equations in Theorem \ref{thm:curvature tensor}.
Also, since this is where our case differs from that in \cite{Belegradek complex}, we need to come up with a new function to vary from $\oshr/2$ to $\cosh(r)$.  
Developing this function (which will just be a cubic polynomial) and proving that we can use Lemma \ref{lemma:general curvature formula} and Corollary \ref{cor:final step} will occupy most of the remainder of this Section and Section \ref{Section:proof of negative curvature}.  

\begin{remark}\label{remark:endpoints}
Note that none of the regions above contain either of its endpoints.  
We deal with this situation as follows.  
For Regions $2$, $4$, and $6$ we will construct smooth functions who, when concatenated with the functions on the surrounding regions, yield $C^1$ functions.  
We will then use Lemma \ref{lemma:warping lemma} (stated below and proved in Appendix A of \cite{Belegradek complex}) to smooth these functions in an arbitrarily small neighborhood of the endpoints. 

With the exceptions of $R_{1414}$, $R_{2424}$, and $R_{3434}$, all of the other components of the curvature tensor depend only on the functions and their first derivatives.  
So if we can control these three components, then we can choose $\d$ small in Lemma \ref{lemma:warping lemma} and apply the Continuity Lemma \ref{lemma:continuity lemma} at these endpoints.
As can be seen in both equation \eqref{eqn:curvature formula} and Lemma \ref{lemma:general curvature formula}, increasing any of $v''$, $h_\th''$, or $h_r''$ decreases the sectional curvature.  
So when bounding the curvature in each region, we will use the smallest values for the respective second derivatives in that region.

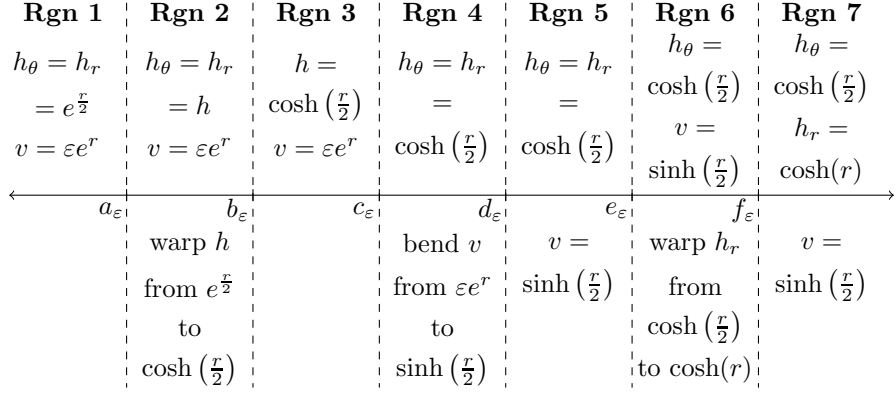
\begin{figure}
\begin{center}
\begin{tikzpicture}[scale=0.8]

\draw[<->] (-1,0) -- (13.7,0);
\draw (0.95,0.1) -- (0.95,-0.1);
\draw (3.05,0.1) -- (3.05,-0.1);
\draw (5.15,0.1) -- (5.15,-0.1);
\draw (7.25,0.1) -- (7.25,-0.1);
\draw (9.35,0.1) -- (9.35,-0.1);
\draw (11.45,0.1) -- (11.45,-0.1);
\draw[dashed] (0.95,3.2) -- (0.95,-3.2);
\draw[dashed] (3.05,3.2) -- (3.05,-3.2);
\draw[dashed] (5.15,3.2) -- (5.15,-3.2);
\draw[dashed] (7.25,3.2) -- (7.25,-3.2);
\draw[dashed] (9.35,3.2) -- (9.35,-3.2);
\draw[dashed] (11.45,3.2) -- (11.45,-3.2);
\draw (0.7,-0.25)node{$a_\e$};
\draw (2.8,-0.25)node{$b_\e$};
\draw (4.9,-0.25)node{$c_\e$};
\draw (7,-0.25)node{$d_\e$};
\draw (9.1,-0.25)node{$e_\e$};
\draw (11.2,-0.25)node{$f_\e$};

\draw (-0.1,3)node{{\bf Rgn 1}};
\draw (2,3)node{{\bf Rgn 2}};
\draw (4.1,3)node{{\bf Rgn 3}};
\draw (6.2,3)node{{\bf Rgn 4}};
\draw (8.3,3)node{{\bf Rgn 5}};
\draw (10.4,3)node{{\bf Rgn 6}};
\draw (12.5,3)node{{\bf Rgn 7}};

\draw (-0.2,2.2)node{$h_\theta = h_r$};
\draw (-0.1,1.5)node{$= e^{\frac{r}{2}}$};
\draw (-0.2,0.8)node{$v = \e e^r$};

\draw (2,2.2)node{$h_\theta = h_r$};
\draw (2,1.5)node{$= h$};
\draw (2,0.8)node{$v = \e e^r$};
\draw (2,-0.8)node{warp $h$};
\draw (2,-1.5)node{from $e^{\frac{r}{2}}$};
\draw (2,-2.2)node{to};
\draw (2,-2.9)node{$\cosh \left( \frac{r}{2} \right)$};

\draw (4.1,2.2)node{$h =$};
\draw (4.1,1.5)node{$\cosh \left( \frac{r}{2} \right)$};
\draw (4.1,0.8)node{$v = \e e^r$};

\draw (6.2,2.2)node{$h_\theta = h_r$};
\draw (6.2,1.5)node{$=$};
\draw (6.2,0.8)node{$\cosh \left( \frac{r}{2} \right)$};
\draw (6.2,-0.8)node{bend $v$};
\draw (6.2,-1.5)node{from $\e e^r$};
\draw (6.2,-2.2)node{to};
\draw (6.2,-2.9)node{$\sinh \left( \frac{r}{2} \right)$};

\draw (8.3,2.2)node{$h_\theta = h_r$};
\draw (8.3,1.5)node{$=$};
\draw (8.3,0.8)node{$\cosh \left( \frac{r}{2} \right)$};
\draw (8.3,-0.8)node{$v =$};
\draw (8.3,-1.5)node{$\inhr/2$};

\draw (10.4,2.5)node{$h_\th =$};
\draw (10.4,1.8)node{$\cosh \left( \frac{r}{2} \right)$};
\draw (10.4,1.1)node{$v =$};
\draw (10.4,0.4)node{$\sinh \left( \frac{r}{2} \right)$};
\draw (10.4,-0.8)node{warp $h_r$};
\draw (10.4,-1.5)node{from};
\draw (10.4,-2.2)node{$\oshr/2$};
\draw (10.4,-2.9)node{to $\cosh (r)$};

\draw (12.5,2.5)node{$h_\theta =$};
\draw (12.5,1.8)node{$\oshr/2$};
\draw (12.5,1.1)node{$h_r =$};
\draw (12.5,0.4)node{$\cosh(r)$};
\draw (12.5,-0.8)node{$v =$};
\draw (12.5,-1.5)node{$\inhr/2$};

\end{tikzpicture}
\end{center}
\caption{Schematic diagram for the warping functions for $\l$, where ``Rgn" stands for ``Region".  Note that the diagram is not remotely close to scale.  The interval $(\ae, \be)$ is very large, while all of the other intervals are arbitrarily small.}
\label{figure:warping stages}
\end{figure}

Lemma \ref{lemma:warping lemma} referenced above is as follows.

\vskip 15pt

\begin{lemma}\label{lemma:warping lemma}
Given real numbers $k, a_1, c, a_2$ with $a_1 < c < a_2$, let $f_1: [a_1,c] \to \R$ and $f_2: [c,a_2] \to \R$ be $C^2$ functions satisfying $f_i'' \geq k$, $f_1(c) = f_2(c)$, and $f_1'(c) \leq f_2'(c)$.  
If $f:[a_1,a_2] \to \R$ denotes the concatenation of $f_1$ and $f_2$, then for any small $\d > 0$ there exists a $C^2$ function $f_\d : [a_1,a_2] \to \R$ such that
	\begin{enumerate}
	\item  $\ds{f_\d'' > k}$.
	\item  $f_\d = f$ and $f_\d' = f'$ at the points $a_1$ and $a_2$.  
	\item  if $f$ is increasing, then $f_\d' > 0$.
	\item  if $f$ is $C^l$ on $[a_1,a_2]$ for some integer $l \in [0,\infty]$, then $f_\d$ if $C^l$ on $[a_1,a_2]$, and $f_\d$ converges to $f$ in the $C^l$-topology on $[a_1,a_2]$ as $\d \to 0$.
	\end{enumerate}
\end{lemma}

\vskip 10pt

Lastly, one added bonus to the method that we are applying at the endpoints is that it greatly simplifies the exposition of Section \ref{Section:proof of negative curvature}.
\end{remark}

\vskip 20pt

\subsection{Simultaneously warping $h_\th$ and $h_r$ from $e^{\frac{r}{2}}$ to $\oshr/2$}\label{subsect:warping h}
Both this and the following Subsection are nearly identical to those in Section 10 of \cite{Belegradek complex}.
The only minor changes are due to the fact that the component of the complex hyperbolic metric $g$ with respect to $\ddt$ is $\inhr/2$ instead of $\sinh(r)$.  

Let $\re$ denote the unique solution to the equation $\e e^r = \inhr/2$.  
Then one sees that $\re \approx 2\e$.  
Let $\re^- = \re - \e^4$, and define $\be = \frac{\re^-}{2}$.  
So notice that $\be = \e + O(\e^3)$ (and please see Remark \ref{rmk:O} for how we will use the ``$O()$" notation).  
The tangent line to the graph of $\oshr/2$ at $\be$ is
	\begin{equation*}
	\ell(r) = \cosh \left( \frac{\be}{2} \right) + \frac{1}{2} \sinh \left( \frac{\be}{2} \right) (r - \be).
	\end{equation*}
Let $q(r) = \ell(r) + \e^6 (r - \be)^2$, and notice that $q(\be) = \ell(\be) = \cosh \left( \frac{\be}{2} \right)$.  

\vskip 15pt

\begin{proposition}\label{prop:warping h}
There exists a $C^1$ function h and values $\ae < m_\e < \be$ such that
	\begin{enumerate}
	\item  The function h is positive and increasing.
	\item  $h(r) = \oshr/2$ for $r \geq \be$.
	\item  $h(r) = q(r)$ for $r \in [m_\e , \be]$.
	\item  If $r \in [\ae, m_\e]$, then h is smooth, $h'' > h/4$, and $(\ln h)'' > 0$ with $(\ln h)' = \frac{h'}{h} \in \left[ \frac{1}{2}, \frac{3}{4} \right]$.
	\item  If $r \leq \ae$, then $h(r) = e^{\frac{r}{2}}$.  
	\end{enumerate}
\end{proposition}

\vskip 10pt

Proposition \ref{prop:warping h} is identical to Proposition 10.4 of \cite{Belegradek complex}, and therefore its proof is omitted here.
One quick remark is that the value of $\ae$ is approximately $- \frac{4}{\e}$.  

\vskip 20pt

\subsection{Bending $v$ from $\e e^r$ to $\inhr/2$}\label{subsect:warping v}
First, recall the definitions of $\re$ and $\re^-$ from the beginning of the previous Subsection.
Let $\ce = \re^-$.
Since the functions $\e e^r$ and $\inhr/2$ intersect at $\re$, their concatenation yields a $C^0$ function.  
The following Proposition just says that we can approximate this by a $C^1$ function whose first two derivatives are controlled nicely.

\vskip 15pt

\begin{proposition}\label{prop:warping v}
There exists a $C^1$ function v and $\de \in (\re, \re + \e^4]$ such that
	\begin{enumerate}
	\item  The function v is positive and increasing.
	\item  $v(r) = \inhr/2$ for $r \geq \de$.
	\item  $v(r) = \e e^r$ for $r \leq \ce$.
	\item  If $r \in [\ce, \de]$, then v is smooth, $v'' > v$, and $(\ln v)'' > 0$.
	\end{enumerate}
\end{proposition}

\vskip 10pt

\begin{remark}
Proposition \ref{prop:warping v} is equivalent to Proposition 10.1 of \cite{Belegradek complex}.
But our situation is slightly different.  
We are bending $\e e^r$ to the function $\inhr/2$, not $\sinh(r)$.  
But one can check that the same proof works with the obvious modifications for our situation.
Every derivative of either $\inhr/2$ or $\oshr/2$ inserts an extra $``\frac{1}{2}"$, but this does not change any of the inequalities within the proof.

So in particular, we still have that ${\bf v}'' > {\bf v}$ on the interval $(\ce,\de)$, even though this inequality is not true for $\inhr/2$.
The function ${\bf v}$ in Propostion \ref{prop:warping v} is only $C^1$, and this shows that it is certainly not $C^2$ at $\de$.
\end{remark}

\vskip 20pt

\subsection{Warping $h_r$ from $\oshr/2$ to $\cosh(r)$}\label{subsect:warping h_r}
Let $\ee = 2 \de$, and let $\fe = (k+1) \ee$ where $k$ is a large positive constant to be chosen independent of $\e$.  
Specifically, we first choose $k$ large and then choose $\e$ small, so that $\fe=(k+1)\ee \to 0$ as $\e \to 0$.  

Define $\v: [\ee, \fe] \to \R$ by
	\begin{equation}\label{eqn:def of varphi}
	\v(r) = C_3 (r - \ee)^3 + C_2 (r-\ee)^2 + \frac{1}{2} \sinh \left( \frac{\ee}{2} \right) (r - \ee) + \cosh \left( \frac{\ee}{2} \right)
	\end{equation}
\vskip 10pt
\n where
	\begin{align*}
	&C_3 = \frac{1}{\d^2} \left( \sinh (\fe) + \frac{1}{2} \sinh \left( \frac{\ee}{2} \right) \right) + \frac{1}{\d^3} \left( 2 \cosh \left( \frac{\ee}{2} \right) - 2 \cosh (\fe) \right)  \\
	&C_2 = - \frac{1}{\d} \left( \sinh(\fe) + \sinh \left( \frac{\ee}{2} \right) \right) + \frac{3}{\d^2} \left( \cosh(\fe) - \cosh \left( \frac{\ee}{2} \right) \right)  \\
	&\d = \fe - \ee = (k+1)\ee - \ee = k \ee.
	\end{align*}
	
\vskip 15pt
	
\n One can check that
	\begin{multicols}{2}
	\n \begin{equation} \v(\ee) = \cosh \left( \frac{\ee}{2} \right)  \label{eqn:varphi 1} \end{equation}
	\begin{equation} \v'(\ee) = \frac{1}{2} \sinh \left( \frac{\ee}{2} \right) \label{eqn:varphi 2}  \end{equation}
	\begin{equation} \v(\fe) = \cosh(\fe)  \label{eqn:varphi 3}  \end{equation}
	\begin{equation} \v'(\fe) = \sinh(\fe)  \label{eqn:varphi 4}  \end{equation}
	\end{multicols}
\n and so one sees that $\v$ is just the cubic polynomial which gives a $C^1$ interpolation between $\oshr/2$ and $\cosh(r)$ on the interval $[\ee, \fe]$.  

The definitions for the constants $C_2$ and $C_3$ are reasonably complicated and can be difficult to use (especially when combined with the formulas in Theorem \ref{thm:curvature tensor}).  
But since $\ee \to 0$ as $\e \to 0$, we can use the Taylor series for $\sinh(r)$ and $\cosh(r)$ about $0$ to give much simpler approximations for these constants.  
But first some notation.

\begin{remark}\label{rmk:O}
We will use the ``big $O$" notation as follows.  
A term is of the order $O(\b)$ if, when divided by $\b$, this term still approaches $0$ as $\b \to 0$.  
We will generally use notation like ``$O(\ee^2)$", which just means that the remaining terms all contain degree three powers (or higher) of $\ee$.
Note that $O(1)$ just means that the remaining terms all contain at least one $\ee$.
We reserve $``\approx"$ for approximations when taking $k$ large.
\end{remark}

With this notation, our estimates for $C_2$ and $C_3$ are as follows:
	\begin{align}
	C_2 &= - \frac{1}{k \ee} \left( \fe + \frac{\ee}{2} + O(\ee^2) \right) + \frac{3}{k^2 \ee^2} \left( 1 + \frac{\fe^2}{2} - 1 - \frac{\ee^2}{8} + O(\ee^3) \right)  \nonumber  \\
	&= \frac{-1}{k \ee} \left( \frac{(2k+3) \ee}{2} \right) + \frac{3}{k^2 \ee^2} \left( \frac{(4k^2 + 8k +3) \ee^2}{8} \right) + O(\ee)  \nonumber \\
	&= \frac{-4k(2k+3) + 12k^2 + 24k + 9}{8k^2} + O(\ee) = \frac{4k^2 + 12k + 9}{8k^2} + O(\ee)  \nonumber 
	\end{align}
	\begin{align}
	&\approx \frac{1}{2} + O(\ee) \qquad \text{ (for } k \text{ large)}  \label{eqn:C2 approx}
	\end{align}
and
	\begin{align}
	C_3 &\approx \frac{1}{k^2 \ee^2} \left( \fe + \frac{\ee}{4} + O(\ee^2) \right) + \frac{1}{k^3 \ee^3} \left( 2 + \frac{\ee^2}{4} - 2 - \fe^2 +O(\ee^3) \right)  \nonumber \\
	&= \frac{1}{k^2 \ee^2} \left( \frac{(4k+5) \ee}{4} \right) + \frac{1}{k^3 \ee^3} \left( \frac{(-4k^2 -8k -3) \ee^2}{4} \right) + O(1)  \nonumber \\
	&= \frac{-3k - 3}{4k^3 \ee} + O(1) = \frac{-3(k+1)}{4k^3 \ee} + O(1)  \nonumber \\
	&\approx \frac{-3}{4 \ee k^2} + O(1) \qquad \text{ (for } k \text{ large)}.  \label{eqn:C3 approx}
	\end{align}
	

	
Now, for $r \in [\ee, \fe]$, we can write $r = C \ee$ for some $1 \leq C \leq k+1$.  
Then letting $\a = C-1$ we have that $r - \ee = \ee (C-1) = \a \ee$ where, of course, $0 \leq \a \leq k$.  
Using this notation with the approximations in equations \eqref{eqn:C2 approx} and \eqref{eqn:C3 approx} turns equation \eqref{eqn:def of varphi} into
	\begin{align}
	\v \approx &\left( \frac{-3}{4 \ee k^2} + O(1) \right) \left( \a^3 \ee^3 \right) + \left( \frac{1}{2} + O(\ee) \right) \left( \a^2 \ee^2 \right)  \nonumber\\
	&+ \frac{1}{2} \left( \frac{\ee}{2} + O(\ee^2) \right) \left( \a \ee) \right) + \left( 1 + \frac{\ee^2}{8} + O(\ee^3) \right)  \nonumber  \\
	= 1 + \ee^2 & \left( \frac{-3}{4k^2} \a^3 + \frac{1}{2} \a^2 + \frac{1}{4} \a + \frac{1}{8} \right) + O(\ee^3).    \label{eqn:approx of varphi}
	\end{align}
	
In Section \ref{Section:proof of negative curvature} we will also need approximations for $\v'$ and $\v''$.  
So let us compute them now.
	\begin{align}
	\v' &= 3 C_3 (r - \ee)^2 + 2C_2 (r - \ee) + \frac{1}{2} \sinh \left( \frac{\ee}{2} \right)  \nonumber  \\
	&\approx  \left( \frac{-9}{4 \ee k^2} + O(1) \right) \left( \a^2 \ee^2 \right) + (1 + O(\ee))(\a \ee) + \left( \frac{\ee}{4} + O(\ee^2) \right)  \nonumber  \\
	&= \ee \left( \frac{-9}{4k^2} \a^2 + \a + \frac{1}{4} \right) + O(\ee^2)  \label{eqn:approx of varphi'}
	\end{align}
and
	\begin{align}
	\v'' &= 6 C_3 (r - \ee) + 2 C_2  \nonumber  \\
	&\approx \left( \frac{-9}{2 \ee k^2} + O(1) \right) (\a \ee) + (1 + O(\ee))  \nonumber \\
	&= \frac{-9}{2k^2} \a + 1 + O(\ee) \approx 1 + O(\ee)  \label{eqn:approx of varphi''}
	\end{align}
since $0 \leq \a \leq k$ and $k$ is large.

\vskip 20pt

\section{Proving the metric in Section \ref{Section:construction of metric} has sectional curvature bounded above by a negative constant}\label{Section:proof of negative curvature}

In this Section we do exactly as the title says, proving Theorem \ref{thm:metric of negative curvature}(1).  
We break our argument up into the different regions.  
Note that region $7$ is clear since the metric $\l$ agrees with the complex hyperbolic metric $g$.
Throughout this Section we again use the notation $R_{ijkl} = \la R(Y_i,Y_j)Y_k,Y_l \ra_\l$.

\vskip 20pt
\subsection*{Region 1}  
Here, $v = \e e^r$ and $h_\th = h_r = e^{\frac{r}{2}} :=h$ over the region $(- \infty, \ae)$.  
Plugging these values into the equations in Corollary \ref{cor:curvature tensor} yields
	\begin{align*}
	&R_{1212} = R_{1313} = - \frac{1}{2} + \frac{\e}{16}  \\
	&R_{2323} = - \frac{1}{4} - \frac{1}{4e^r} - \frac{3\e^2}{16}  \\
	&R_{1414} = - 1 \hskip 80pt R_{2424} = R_{3434} = - \frac{1}{4}  \\
	&R_{1423} = \frac{-\e}{4}   \hskip 80pt  R_{1234} = \frac{\e}{8} = - R_{1324}.
	\end{align*}
Then Lemma \ref{lemma:general curvature formula} proves that the sectional curvature of $\l$ is bounded above by a negative constant on this region.

\vskip 20pt
\subsection*{Region 2}
In this region $v = \e e^r$ and $h_\th = h_r :=h$ is warped from $e^{\frac{r}{2}}$ to $\oshr/2$ over the interval $(\ae, \be)$.  
Region 2 is broken up into two different situations:  $r \in (\ae, m_\e)$ and $r \in (m_\e, \be)$.  
We deal with the two intervals separately.  
Our arguments for these regions are virtually identical to those in \cite{Belegradek complex}.  
We include the arguments here because our curvature formulas in Corollary \ref{cor:curvature tensor} are slightly different than Belegradek's, and to verify a few estimates in our case (since our $r_\e$ is twice that of what is in \cite{Belegradek complex}).  

\subsection*{The interval $(\ae, m_\e)$}
By Proposition \ref{prop:warping h} we know that $\frac{h''}{h} > \frac{1}{4}$, $\frac{h'}{h}$ is increasing (since $(\ln h)'' > 0$), and $\frac{h'}{h} \in [ \frac{1}{2}, \frac{3}{4} ]$.  
Of course $\frac{v''}{v} = \frac{v'}{v} = 1$, and so $\frac{v'}{v} - \frac{h'}{h} \leq \frac{1}{2} < 1$.  
Since $\ln h$ is convex, its graph lies above its tangent line at $\ae$.  
Hence $\ln h \geq \frac{r}{2}$ implying that $h \geq e^{\frac{r}{2}}$.  
Therefore, $\frac{v}{h^2} \leq \frac{\e e^r}{e^r} = \e < 2 \e$.  

Plugging these estimates into the formulas in Corollary \ref{cor:curvature tensor} yields
	\begin{align*}
	&K(Y_1,Y_2) = K(Y_1,Y_3) \leq - \frac{1}{2} + \frac{\e^2}{4} < - \frac{1}{3}  \\
	&K(Y_2,Y_3) < - \left( \frac{h'}{h} \right)^2 \leq - \frac{1}{4}  \\
	&K(Y_1,Y_4) = - 1 \hskip 70pt K(Y_2,Y_4) = K(Y_3,Y_4) = - \frac{h''}{h} < - \frac{1}{4}  \\
	&|R_{1423}| = \biggr| \frac{-v}{2h^2} \left( \frac{v'}{v} - \frac{h'}{h} \right) \biggr| < \e  \hskip 50pt  |R_{1234}| = |R_{1324}| < \frac{\e}{2}.
	\end{align*}
Lemma \ref{lemma:general curvature formula} then completes the argument for $(\ae, m_\e)$.

\subsection*{The interval $(m_\e, \be)$}
Over this interval $v = \e e^r$ and $h_\th = h_r = q(r) := h$, where $q(r)$ is defined in Subsection \ref{subsect:warping h}.  
Two things that are immediately clear are that $\frac{v''}{v} = 1$ and $\frac{h''}{h} > \e^6$.  
We also have that
	\begin{equation*}
	q'(r) = \frac{1}{2} \sinh \left( \frac{\be}{2} \right) + 2 \e^6 (r - \be) = \frac{\e}{4} + O(\e^2)
	\end{equation*}
and so, in particular, $q$ is increasing.  
Thus $q(r) < q(\be) = \cosh \left( \frac{\be}{2} \right) = 1 + O(\e^2)$.  
Since $q'' = 2\e^6$, one easily checks that $\frac{q'}{q}$ is decreasing over the interval $(m_\e, \be)$ from $\frac{3}{4}$ to $\frac{1}{2} \tanh \left( \frac{\be}{2} \right) = \frac{\e}{4} + O(\e^2)$.  
Therefore, $\frac{h'}{h} \in \left( \frac{\e}{5}, \frac{4}{5} \right)$.  

The last quantity that we need to bound is $\frac{v}{h^2} = \frac{\e e^r}{q^2}$.  
The argument is identical to pg. 567 of \cite{Belegradek complex}, and so we omit it here.  
The idea of the argument is that you differentiate $\frac{v}{h^2}$ twice to show that it is locally maximized at the endpoints $\{ m_\e, \be \}$, and then argue that the maximum value is actually at $r = \be$.  
You then have that
	\begin{equation*}
	\frac{v}{h^2} \leq \e \, \frac{e^{\be}}{\cosh^2 \left( \frac{\be}{2} \right)} < 2 \e.
	\end{equation*}
Plugging these estimates into the equations in Corollary \ref{cor:curvature tensor} gives
	\begin{align*}
	&K(Y_1,Y_2) = K(Y_1,Y_3) < - \frac{\e}{5} + \frac{\e^2}{4} < - \frac{\e}{6}   \\
	&K(Y_2,Y_3) < - \frac{1}{4h^2} \leq \frac{-1}{4 \cosh^2 \left( \frac{\be}{2} \right)} < - \frac{1}{9}  \\
	&K(Y_1,Y_4) = -1  \hskip 80pt K(Y_2,Y_4) = K(Y_3,Y_4) < - \e^6  \\
	&|R_{1423} = \biggr| \frac{-v}{2h^2} \left( \frac{v'}{v} - \frac{h'}{h} \right) \biggr| < \frac{1}{2} \cdot (2 \e) \left( 1 - \frac{\e}{5} \right) < \e.
	\end{align*}
Then inserting these values into equation \eqref{eqn:curvature formula} yields
	\begin{align*}
	K(\s) < - \frac{\e}{6} &\left((c_1d_2 - c_2d_1)^2 + d_1^2c_3^2 \right) - d_1^2c_4^2  - \frac{1}{9} d_2^2c_3^2 - \e^6 d_2^2c_4^2 + 3 \e |c_3c_4d_1d_2|  \\
	= - \frac{\e}{6} &\left((c_1d_2 - c_2d_1)^2 + d_1^2c_3^2 \right) - \e^6 d_2^2c_4^2 - \left( |c_4 d_1| - \frac{1}{3} |c_3 d_2| \right)^2  \\
	&- \left( \frac{2}{3} - 3 \e \right) | c_3 c_4 d_1 d_2 | .
	\end{align*}
Every term in the above sum is nonpositive, and not all of the coefficients can simultaneously be zero.  
So by compactness we have that $K()$ is bounded above by a negative constant within Region $(m_\e, \be)$.

\vskip 20pt
\subsection*{Region 3}
Here, $v = \e e^r$ and $h_\th = h_r = \oshr/2:=h$ over the region $( \be, \ce)$.  
Recall that $\be \approx \e$ and $\ce \approx 2\e$.
Then for $\e > 0$ small enough we have:
	\begin{itemize}
	\item  $\ds{\frac{v}{h^2} = \frac{\e e^r}{\cosh^2 \left( \frac{r}{2} \right)} \leq \frac{\e e^{2 \e}}{1} < 2 \e.}$
	\vskip 6pt
	\item  $\ds{\frac{h'}{h} = \frac{\inhr/2}{2 \oshr/2} \geq \frac{\frac{r}{2}}{4} = \frac{r}{8} \geq \frac{\e}{8} > \frac{\e}{9}.}  $
	\end{itemize}
	\vskip 6pt
	
Then we can plug into the equations in Corollary \ref{cor:curvature tensor} to obtain:
	\begin{align*}
	&K(Y_1,Y_2) = K(Y_1,Y_3) = - \frac{h'}{h} + \frac{v^2}{16 h^4} < \frac{\e}{9} + \frac{\e^2}{4} < - \frac{\e}{10}  \\
	&K(Y_2,Y_3) < \frac{-1}{4h^2} < - \frac{1}{9}  \\
	&K(Y_1,Y_4) = -1  \hskip 80pt K(Y_2,Y_4) = K(Y_3,Y_4) = - \frac{1}{4}  \\
	&|R_{1423}| = \frac{1}{2} \frac{v}{h^2} \left( \ln \frac{v}{h} \right)^\prime \leq \frac{1}{2} (2 \e) \left( \frac{v'}{v} - \frac{h'}{h} \right) = \frac{3}{4} \e < \e.
	\end{align*}
	
\vskip 8pt
	
\n To find $K(\s)$ for any $2$-plane $\s$ we plug these values into equation \eqref{eqn:curvature formula}, giving
	\begin{align*}
	K(\s) \leq - \frac{\e}{10} & \left( (c_1 d_2 - c_2 d_1)^2 + d_1^2 c_3^2 \right) - \frac{1}{4} (d_1^2 c_4^2 + d_2^2 c_4^2) - \frac{1}{9} d_2^2 c_3^2 + 3\e |c_3 c_4 d_1 d_2|  \\
	= - \frac{\e}{10} & \left( (c_1 d_2 - c_2 d_1)^2 + d_1^2 c_3^2 \right) - \frac{1}{4} d_2^2 c_4^2 - \left( \frac{1}{2} |c_4 d_1| - \frac{1}{3} |c_3 d_2| \right)^2  \\
	&+ |c_3c_4d_1d_2|\left( 3 \e - \frac{1}{3} \right).
	\end{align*}
Every term in the sum above is nonpositive.  
So we have that $K()$ is bounded above by a negative constant within Region $3$.

\vskip 20pt
\subsection*{Region 4}
In this region $h_\th = h_r = \oshr/2$ over the region $(\ce, \de)$, while $v$ is ``smoothed out" from $\e e^r$ to $\inhr/2$.  
Recall from Proposition \ref{prop:warping v} that over this region $v$ and $v'$ are increasing, $v'' \geq \frac{v}{4}$, and $(\ln v)'' > 0$ (which, in particular, implies that $\frac{v'}{v}$ is increasing).
Also, recall that $\ce = \re^- = \re - \e^4 \approx 2 \e - \e^4$ and $\de \leq \re + \e^4 \approx 2 \e + \e^4 $.
So over this entire region, $r \approx 2 \e$.

Since both $v$ and $h$ are increasing we have that
	\begin{equation*}
	\frac{v}{h^2} \leq \frac{v(\de)}{h(\ce)} \approx \frac{\sinh \left( \frac{2 \e + \e^4}{2} \right)}{\cosh \left( \frac{2 \e - \e^4}{2} \right)} \approx \frac{\frac{2 \e}{2}}{1} = \e < 2 \e.
	\end{equation*}

We also know that $\frac{v'}{v}$ is increasing, and so it can be bounded by its values at the endpoints of the interval (where $v = \e e^r$ and $\inhr/2$, respectively).  
Therefore,
	\begin{equation*}
	1 \leq \frac{v'}{v} \leq \frac{1}{2} \coth \left( \frac{\de}{2} \right)
	\end{equation*}
which, for $\e > 0$ small, is very large.  
We also have that
	\begin{equation*}
	\frac{h'}{h} = \frac{\inhr/2}{2 \oshr/2} = \frac{1}{2} \tanh \left( \frac{r}{2} \right) \approx \frac{r}{4} \approx \frac{2 \e}{4} = \frac{\e}{2}.
	\end{equation*}
Combining these last two equations gives us
	\begin{itemize}
	\item  $\ds{ 0 < \frac{v'}{v} - \frac{h'}{h} < \frac{1}{2} \coth \left( \frac{\de}{2} \right) }$
	\vskip 6pt
	\item  $\ds{ \frac{h' v'}{hv} > \frac{\e}{4}  }$
	\vskip 6pt
	\item  $\ds{ \frac{v}{h^2} \left( \frac{v'}{v} - \frac{h'}{h} \right) \leq \frac{\sinh \left( \frac{\de}{2} \right)}{\cosh \left( \frac{\ce}{2} \right) } \cdot \frac{1}{2} \coth \left( \frac{\de}{2} \right) = \frac{\cosh \left( \frac{\de}{2} \right)}{2 \cosh \left( \frac{\ce}{2} \right)} = \frac{1}{2} + O(\e^2). }$
	\end{itemize}
	
\vskip 10pt
	 
\n Combining these estimates with the formulas in Corollary \ref{cor:curvature tensor} gives us
	\begin{align*}
	&K(Y_1,Y_2) = K(Y_1,Y_3) = - \frac{v' h'}{vh} + \frac{v^2}{16h^4} < -\frac{\e}{4} + \frac{\e^2}{4} < - \frac{\e}{5}  \\
	&K(Y_2,Y_3) < - \frac{1}{4h^2} - \left( \frac{h'}{h} \right)^2 = \frac{-1 - \sinh^2 \left( \frac{r}{2} \right)}{4 \cosh^2 \left( \frac{r}{2} \right)} = - \frac{1}{4}  \\
	&K(Y_1,Y_4) < - 1 \hskip 80pt K(Y_2,Y_4) = K(Y_3,Y_4) = - \frac{1}{4}  \\
	&|R_{1423}| = \biggr| \frac{v}{2h^2} \left( \frac{v'}{v} - \frac{h'}{h} \right) \biggr|  < \frac{1}{4} + \e.
	\end{align*}
	
\vskip 8pt
	
\n To find $K(\s)$ for any $2$-plane $\s$ we plug these values into equation \eqref{eqn:curvature formula}, giving
	\begin{align*}
	K(\s) &< - \frac{\e}{5} \left( (c_1d_2 - c_2d_1)^2 + c_3^2d_1^2 \right) - c_4^2 d_1^2 - \frac{1}{4} (c_3^2d_2^2+ c_4^2d_2^2) + 3 \left( \frac{1}{4} + \e \right) |c_3c_4d_1d_2|  \\
	&= - \frac{\e}{5} \left( (c_1d_2 - c_2d_1)^2 + d_1^2c_3^2 \right) - \left( c_4d_1 - \frac{1}{2} d_2 c_3 \right)^2 - \frac{1}{4} c_4^2d_2^2 \\
	&- \left(\frac{5}{4} - \e \right) |c_3c_4d_1d_2|.
	\end{align*}
Every term in the sum above is nonpositive, and so we have that $K()$ is bounded above by a negative constant within Region $4$.

\vskip 20pt
\subsection*{Region 5}
Here, $v = \inhr/2$ and $h_\th = h_r = \oshr/2 :=h$ over the region $(\de, \ee)$.  
Plugging these values into the equations in Corollary \ref{cor:curvature tensor} yields
	\begin{align*}
	&R_{1212} = R_{1313} = - \frac{1}{4} + \frac{\sinh^2 \left( \frac{r}{2} \right)}{16 \cosh^4 \left( \frac{r}{2} \right)} \approx - \frac{1}{4}  \\
	&R_{2323} = - \frac{1}{4} - \frac{3 \sinh^2 \left( \frac{r}{2} \right)}{16 \cosh^4 \left( \frac{r}{2} \right)} < - \frac{1}{4}  \\
	&R_{1414} = R_{2424} = R_{3434} = - \frac{1}{4}  \\
	&R_{1423} = \frac{- 1}{4 \cosh^3 \left( \frac{r}{2} \right)} > - \frac{1}{4}   \hskip 80pt  R_{1234} = -R_{1324} = \frac{1}{8 \cosh^3 \left( \frac{r}{2} \right)} \approx \frac{1}{8}.
	\end{align*}
Then for $\e > 0$ small enough, Lemma \ref{lemma:general curvature formula} proves that the sectional curvature of $\l$ is bounded above by a negative constant on this region.

\vskip 20pt
\subsection*{Region 6}
In this region $v = \inhr/2$, $h_\th = \oshr/2$, and $h_r = \v$ where $\v$ is the cubic defined in Subsection \ref{subsect:warping h_r} which varies in a $C^1$ manner between $\oshr/2$ and $\cosh(r)$.  
If one plugs in the values $v = \inhr/2$, $h_\th = \oshr/2$ and $h_r = \v$ into the equations in Theorem \ref{thm:curvature tensor} and then simplifies (a lot), they come up with the following.
	\begin{align*}
	&R_{1212} = \frac{- (\v^2-1)(3\v^2 + \cosh^2(r))}{4 \v^2 \sinh^2(r)}  \\
	&R_{1313} = \frac{- (\cosh(r) + 1) \v'}{2 \sinh(r) \v} - \frac{1}{4 \v^2} + \frac{(\v^2-1)(\v^2 + 2 \cosh(r) + 1)}{4 \v^2 \sinh^2(r)}  \\
	&R_{2323} = \frac{- (\cosh(r) - 1) \v'}{2 \sinh(r) \v} - \frac{1}{4 \v^2} + \frac{(\v^2-1)(\v^2 - 2 \cosh(r) + 1)}{4 \v^2 \sinh^2(r)}  \\
	&R_{1414} = - \frac{1}{4} \hskip 40pt R_{2424} = - \frac{1}{4} \hskip 40pt R_{3434} = - \frac{\v''}{\v}  \\
	&R_{1234} = \frac{\v'}{\sinh(r)} - \frac{(\v^2-1) \cosh(r)}{2 \v \sinh^2(r)}  \\
	&R_{1324} = \frac{\v' (\cosh(r) + \v^2)}{2 \v^2 \sinh(r)} - \frac{(\v^2-1)(\cosh(r) + 1)}{4 \v \sinh^2(r)} - \frac{1}{4 \v}  \\
	&R_{1423} =  \frac{\v' (\cosh(r) - \v^2)}{2 \v^2 \sinh(r)} + \frac{(\v^2 - 1)(\cosh(r) - 1)}{4 \v \sinh^2(r)} - \frac{1}{4 \v} .
	\end{align*}
One nice way to ``check" these formulas is to let $\v = \cosh(r)$ and confirm that this gives you the values in equations \eqref{eqn:curvature 1} through \eqref{eqn:curvature 9}.

Our method to prove that $K()$ is bounded above by a negative constant in region $6$ is to attempt to apply Lemma \ref{lemma:general curvature formula}, and when exactly one of these inequalities fails by an arbitrarily small amount to really apply Corollary \ref{cor:final step}.  
So we need to consider the inequalities in Lemma \ref{lemma:general curvature formula}.  
But first we compute two estimates (equations \eqref{eqn:approx of varphi'/sinh} and \eqref{eqn:approx of varphi2 - 1/sinh2}) which show up in many of the inequalities.

The Taylor series for $\sinh(r)$ centered at $r = \ee$ is
	\begin{align*}
	\sinh(r) &= \sinh(\ee) + \cosh(\ee) (r - \ee) + \frac{\sinh(\ee)}{2} (r-\ee)^2 +  \hdots  \\
	&= \left( \ee + \frac{\ee^3}{6} \right) + \left( 1 + \frac{\ee^2}{2} \right) \ee \a + \frac{1}{2} \left( \ee + \frac{\ee^3}{6} \right) \ee^2 \a^2 +  O(\ee^2)  \\
	&= \ee(1 + \a) + O(\ee^2).
	\end{align*}
Combining this with equation \eqref{eqn:approx of varphi'} gives
	\begin{equation}\label{eqn:approx of varphi'/sinh}
	\frac{\v'(r)}{\sinh(r)} \approx \frac{\ee \left( \frac{1}{4} + \a - \frac{9}{4k^2} \a^2 \right) + O(\ee^2)}{\ee (1 + \a) + O(\ee^2)} = \frac{\frac{1}{4} + \a - \frac{9}{4k^2} \a^2}{1 + \a} + O(\ee).
	\end{equation}

\vskip 10pt

\n Now, using equation \eqref{eqn:approx of varphi} one sees that
	\begin{equation*}
	\v^2(r) - 1 \approx \ee^2 \left( \frac{1}{4} + \frac{1}{2} \a + \a^2 - \frac{3}{2k^2} \a^3 \right) + O(\ee^3).
	\end{equation*}
Also, the above formula for $\sinh(r)$ gives us that
	\begin{equation*}
	\sinh^2(r) = \ee^2 (1 + \a)^2 + O(\ee^3).
	\end{equation*}
Therefore
	\begin{equation}\label{eqn:approx of varphi2 - 1/sinh2}
	\frac{\v^2(r) - 1}{\sinh^2(r)} = \frac{\frac{1}{4} + \frac{1}{2} \a + \a^2 - \frac{3}{2k^2} \a^3}{(1+\a)^2} + O(\ee).
	\end{equation}
	
\vskip 10pt

We now deal with each of the six inequalities in Lemma \ref{lemma:general curvature formula}.  
We consider them in reverse order since the inequality that will cause issues is inequality $(1a)$, and we will prove that each inequality (except (1a)) is strict.
Note that we will have to derive estimates other than \eqref{eqn:approx of varphi'/sinh} and \eqref{eqn:approx of varphi2 - 1/sinh2} in some of these situations, but these two come up so often that we did them first.

\vskip 20pt
\subsubsection*{Inequality (3b)}  
Inequality (3b) is $R_{2323} < - |R_{1423}|$.   
Over region $6$ we have that $R_{1423} < 0$, and so we show that $R_{2323} < R_{1423}$.  
Using the equations above, we see that this inequality holds if and only if
	\begin{align*}
	\v - 1 &< 2 \left( \frac{\v'}{\sinh(r)} \right) (1 + \v) (\cosh(r) - \v)  \\
	&- \left( \frac{\v^2 - 1}{\sinh^2(r)} \right) (\v^2 - 2 \cosh(r) + 1 - \v \cosh(r) + \v)
	\end{align*}

\n To prove the above inequality, we need the following estimates:
	\begin{align*}
	\cosh(r) &= 1 + \ee^2 \left( \frac{1}{2} + \a + \frac{1}{2} \a^2 \right) + O(\ee^3)  \\
	\v^2(r) &\approx 1 + \ee^2 \left( \frac{1}{4} + \frac{1}{2} \a + \a^2 - \frac{3}{2k^2} \a^3 \right) + O(\ee^3)  \\
	\cosh(r) - \v(r) &\approx \ee^2 \left( \frac{3}{8} + \frac{3}{4} \a + \frac{3}{4k^2} \a^3 + O(\ee) \right)  \\
	\v^2(r) - 2 \cosh(r) + 1 - \v \cosh(r) + \v &\approx  \ee^2 \left( - \frac{5}{4} - \frac{5}{2} \a - \frac{1}{2} \a^2 - \frac{3}{2k^2} \a^3 + O(\ee) \right).
	\end{align*}
Then, applying equations \eqref{eqn:approx of varphi'/sinh} and \eqref{eqn:approx of varphi2 - 1/sinh2} with the above estimates, we have that inequality (3b) holds if and only if
	\begin{align*}
	\ee^2 &\left( \frac{1}{8} + \frac{1}{4} \a + \frac{1}{2} \a^2 - \frac{3}{4k^2} \a^3 \right) < \ee^2 \left( \frac{\frac{1}{4} + \a - \frac{9}{4k^2}\a^2}{1 + \a} \right) \left( \frac{3}{2} + 3 \a + \frac{3}{k^2} \a^3 \right) \\
	&+ \ee^2 \left( \frac{\frac{1}{4} + \frac{1}{2} \a + \a^2 - \frac{3}{2k^2} \a^3}{(1+\a)^2} \right) \left( \frac{5}{4} + \frac{5}{2} \a + \frac{1}{2} \a^2 + \frac{3}{2k^2} \a^3 \right) + O(\ee^3).
	\end{align*}
This inequality is true if and only if
	\begin{align*}
	O(\ee) < \frac{9}{16} + \frac{27}{8} \a &+ \left( \frac{27}{4} - \frac{9}{4k^2} \right) \a^2 + \left( \frac{9}{2} - \frac{9}{k^2} \right) \a^3 - \frac{3}{4k^2} \a^4  \\
	&+ \left( \frac{18}{4k^2} - \frac{27}{4k^4} \right) \a^5 - \frac{9}{k^4} \a^6  \\
	\iff \hskip 30pt O(\ee) < \frac{9}{16} + \frac{27}{8} \a &+ \left( \frac{27}{4} - \frac{9}{4k^2} \right) \a^2 + \left( \frac{9}{4} - \frac{9}{k^2} \right) \a^3 + \left( \frac{9}{4} - \frac{3\a}{4k^2} \right) \a^3  \\
	&+ \left( \frac{9}{4k^2} - \frac{27}{4k^4} \right) \a^5 + \left( \frac{9}{4k^2} - \frac{9 \a}{k^4} \right)\a^5.
	\end{align*}
For $k$ large enough, every term on the right hand side of the above equation is positive (since $0 \leq \a \leq k$).  
So the right hand side is bounded below by $\frac{9}{16}$.
This proves that inequality (3b) holds for $\e$ sufficiently small.

\vskip 20pt
\subsubsection*{Inequality (3a)}  
Inequality (3a) is $R_{1414} < - |R_{1423}|$.   
Over region $6$ we have that $R_{1423} < 0$, and so we show that $R_{1414} < R_{1423}$.  
Using the equations above, we see that this inequality holds if and only if
	\begin{align*}
	\frac{\cosh(r) - \v^2}{\v(1 - \v)} &< \frac{\sinh(r)}{\v'} \left( \frac{\v + \cosh(r) + 2}{2(\cosh(r) + 1)} \right)  \\
	&= \frac{\sinh(r)}{\v'} + O(\ee).
	\end{align*}

To prove the above inequality we need the following estimates:
	\begin{equation*}
	\cosh(r) - \v^2(r) \approx \ee^2 \left( \frac{1}{4} + \frac{1}{2} \a - \frac{1}{2} \a^2 + \frac{3}{2k^2} \a^3 + O(\ee) \right)
	\end{equation*}
and
	\begin{equation*}
	\v(r) (1 - \v(r)) \approx \ee^2 \left( - \frac{1}{8} - \frac{1}{4} \a - \frac{1}{2} \a^2 + \frac{3}{4k^2} \a^3 + O(\ee) \right).
	\end{equation*}
These can be obtained by using the estimates derived for inequality (3b) along with equation \eqref{eqn:approx of varphi}.
Then, using the above estimates with equation \eqref{eqn:approx of varphi'/sinh}, inequality (3a) is satisfied if and only if
	\begin{align*}
	&\frac{\frac{1}{4} + \frac{1}{2} \a - \frac{1}{2} \a^2 + \frac{3}{2k^2} \a^3}{- \frac{1}{8} - \frac{1}{4} \a - \frac{1}{2} \a^2 + \frac{3}{4k^2} \a^3} < \frac{1 + \a}{\frac{1}{4} + \a - \frac{9}{4k^2} \a^2} + O(\ee)  \\
	\iff \hskip 8pt O(\ee) &< 6 + 24 \a + \left( 36 - \frac{18}{k^2} \right) \a^2 - \frac{48}{k^2} \a^3 + \frac{60}{k^2} \a^4 - \frac{108}{k^4} \a^5  \\
	\iff \hskip 8pt O(\ee) &< 6 + 24 \a + \left( 18 - \frac{18}{k^2} \right) \a^2 + \left( 18 - \frac{48}{k^2} \a \right) \a^2 + \frac{1}{k^2} \left( 60 - \frac{108}{k^2} \a \right) \a^4.
	\end{align*}
Recall that $0 \leq \a \leq k$.  
So every term on the right hand side of the last inequality is positive for $k$ sufficiently large.  
Therefore the right hand side of the inequality is bounded below by $6$, which verifies that inequality (3a) is satisfied for $\e$ sufficiently small.

\vskip 20pt
\subsubsection*{Inequality (2b)}  
Inequality (2b) is $R_{2424} < - |R_{1324}|$.   
Over region $6$ we have that $R_{1324}$ varies from $\approx - \frac{1}{8}$ to $\frac{1}{4}$.
We will show that $R_{2424} < - R_{1324}$, and an analogous (simpler) argument shows that $R_{1414} < R_{1324}$.   
Using the equations above, we see that this inequality holds if and only if
	\begin{align*}
	- \frac{1}{4} &< - \frac{\v' (\cosh(r) + \v^2)}{2 \v^2 \sinh(r)} + \frac{(\v^2-1)(\cosh(r) + 1)}{4 \v \sinh^2(r)} + \frac{1}{4 \v}  \\
	&= - \frac{\v'}{\sinh(r)} + \frac{(\v^2 - 1)}{2 \sinh^2(r)} + \frac{1}{4} + O(\ee).
	\end{align*}

Using equations \eqref{eqn:approx of varphi'/sinh} and \eqref{eqn:approx of varphi2 - 1/sinh2}, we see that inequality (2b) is satisfied if and only if
	\begin{align*}
	-\frac{1}{2} + O(\ee) &< - \left( \frac{\frac{1}{4} + \a - \frac{9}{4k^2} \a^2}{1 + \a} \right) + \frac{1}{2} \left( \frac{\frac{1}{4} + \frac{1}{2} \a + \a^2 - \frac{3}{2k^2} \a^3}{(1+\a)^2} \right)  \\
	\iff \hskip 20pt - (1 + \a)^2 + O(\ee) &< - \frac{1}{4} - 2 \a + \left( \frac{9}{2k^2} - 1 \right) \a^2 + \frac{3}{k^2} \a^3  \\
	\iff \hskip 76pt O(\ee) &< \frac{3}{4} + \frac{9}{2k^2} \a^2 + \frac{3}{k^2} \a^3.
	\end{align*}
The right hand side of this inequality is clearly bounded below by $\frac{3}{4}$ since $0 \leq \a \leq k$.  
Thus, for $\e$ sufficiently small, inequality (2b) is satisfied.

\vskip 20pt
\subsubsection*{Inequality (2a)}  
Inequality (2a) is $R_{1313} < - |R_{1324}|$.   
Over region $6$ we have that $R_{1324}$ varies from $\approx - \frac{1}{8}$ to $\frac{1}{4}$.
We will show that $R_{1313} < - R_{1324}$, and an analogous (simpler) argument shows that $R_{1313} < R_{1324}$.   
Using the equations above, we see that this inequality holds if and only if
	\begin{align*}
	- \frac{(\v + 1)}{4 \v^2} &< \frac{\v' (\v + \v \cosh(r) - \cosh(r) - \v^2)}{2 \v^2 \sinh(r)}  \\
	&+ \frac{(\v^2-1)(\v \cosh(r) + \v - \v^2 - 2 \cosh(r) - 1)}{4 \v^2 \sinh^2(r)}  \\
	\iff \hskip 20pt - \frac{1}{2} + O(\ee) &< \left( \frac{\v'}{2\sinh(r)} \right) (\v + \v \cosh(r) - \cosh(r) - \v^2) - \frac{(\v^2-1)}{2 \sinh^2(r)}.
	\end{align*}
But one can check that
	\begin{equation*}
	\v + \v \cosh(r) - \cosh(r) - \v^2 = O(\ee)
	\end{equation*}
and therefore what we need to show is that
	\begin{equation*}
	\frac{(\v^2-1)}{\sinh^2(r)} + O(\ee) < 1.
	\end{equation*}
Using equation \eqref{eqn:approx of varphi2 - 1/sinh2}, we have that inequality (2a) is satisfied if and only if
	\begin{align*}
	\frac{\frac{1}{4} + \frac{1}{2} \a + \a^2 - \frac{3}{2k^2} \a^3}{(1+\a)^2} + O(\ee) &< 1  \\
	\iff \hskip 30pt O(\ee) &< \frac{3}{4} + \frac{3}{2} \a + \frac{3}{2k^2} \a^3.
	\end{align*}
The right hand side is clearly bounded below by $\frac{3}{4}$ for $0 \leq \a \leq k$, proving that inequality (2a) is satisfied for $\e$ sufficiently small.

\vskip 20pt
\subsubsection*{Inequality (1b)}  
Inequality (1b) is $R_{3434} < - |R_{1234}|$.   
Over region $6$ we have that $R_{1234} > 0$, and so we will show that $R_{3434} < - R_{1234}$.  
Using the equations above, we see that this inequality holds if and only if
	\begin{align*}
	 - \frac{\v''}{\v} &< - \frac{\v'}{\sinh(r)} + \frac{(\v^2-1)\cosh(r)}{2 \v \sinh^2(r)}  \\
	\iff \hskip 30pt -1 + O(\ee) &< - \frac{\v'}{\sinh(r)} + \frac{\v^2-1}{2 \sinh^2(r)} + O(\ee)  \\
	\iff \hskip 30pt -1 + O(\ee) &< - \frac{\v'}{\sinh(r)} + \frac{\v^2-1}{2 \sinh^2(r)}.
	\end{align*}

Note that this inequality is weaker than inequality (2b).  
For (2b), we had a $``- \frac{1}{2}"$ instead of a $``-1"$ on the left hand side of the inequality.  
Therefore, since inequality (2b) holds for $\e$ sufficiently small, so does inequality (1b).

\vskip 20pt
\subsubsection*{Inequality (1a)}  
Inequality (1a) is $R_{1212} < - |R_{1234}|$.   
Over region $6$ we have that $R_{1234} > 0$, and so we will (try to) show that $R_{1212} < - R_{1234}$.  
Using the equations above, we see that this inequality holds if and only if
	\begin{align*}
	\frac{- (\v^2-1)(3\v^2 + \cosh^2(r))}{4 \v^2 \sinh^2(r)} &< - \frac{\v'}{\sinh(r)} + \frac{(\v^2-1) \cosh(r)}{2 \v \sinh^2(r)}  \\
	\iff \hskip 30pt - \frac{(\v^2-1)}{\sinh^2(r)} + O(\ee) &< - \frac{\v'}{\sinh(r)} + \frac{\v^2-1}{2 \sinh^2 (r)} + O(\ee)  \\
	\iff \hskip 43pt \frac{\v'}{\sinh(r)} + O(\ee) &< \frac{3}{2} \left( \frac{\v^2 - 1}{\sinh^2(r)} \right) .
	\end{align*}
Using the above estimates, inequality (1a) is satisfied if and only if
	\begin{align*}
	\frac{\frac{1}{4} + \a - \frac{9}{4k^2} \a^2}{1 + \a} + O(\ee) &< \frac{3}{2} \left( \frac{\frac{1}{4} + \frac{1}{2} \a + \a^2 - \frac{3}{2k^2} \a^3}{(1+\a)^2} \right)  \\
	\iff \hskip 30pt  O(\ee) &< \frac{1}{4} - \a + \left( 1 + \frac{9}{2k^2} \right) \a^2.
	\end{align*}
Let $p(\a) = \frac{1}{4} - \a + \left( 1 + \frac{9}{2k^2} \right) \a^2.$  
One can check that $p$ is an upward opening parabola which is minimized at the point
	\begin{equation*}
	\a^* = \frac{k^2}{2k^2 + 9}.
	\end{equation*}  
Also, one can check that
	\begin{equation*}
	p(\a^*) = \frac{9}{4(2k^2+9)} > 0.
	\end{equation*}
	
\vskip 10pt
	
This seems good, but there is a problem.  
The problem is that $p(\a^*) \to 0$ as $k \to \infty$.  
The reason that this is (potentially) a problem is because we estimated $C_2$ and $C_3$ in equations \eqref{eqn:C2 approx} and \eqref{eqn:C3 approx} for $k$ large, and then used these estimates for $\v$ and $\v'$.  
It is possible that we could have $p(\a^*) < 0$ when we plug in the actual values for $C_2$ and $C_3$.  

We fix this as follows.  
We may consider the components of the curvature tensor $R_{ijkl}$ as functions of the independent variables $\ee$, $\a$, and $k$, where $\ee, k > 0$ and $0 \leq \a \leq k$.  
The first five inequalities are all satisfied for all values of $\a$ if $\e$ is sufficiently small and $k$ is sufficiently large.  
So there exists $\e_1, k_1 > 0$ such that these inequalities hold for all $\e \leq \e_1$ and $k \geq k_1$.  

Notice that $\a^* \to \frac{1}{2}$ as $k \to \infty$.  
Set $\a = \frac{1}{2}$ in the first five inequalities (inequalities (1b) through (3b)) and fix $\e \leq \e_1$.  
For each $k \geq k_1$ let $\d_k > 0$ be the small positive constant by which inequality (1a) is permitted to fail (whose existence is guaranteed by Corollary  \ref{cor:final step}).
Since inequalities (1b) through (3b) all hold as $k \to \infty$, there exists some corresponding $\d_\infty > 0$ also guaranteed by Corollary \ref{cor:final step}.  
Clearly, $\d_k \to \d_\infty$ as $k \to \infty$.  

Since $p(\a^*) \to 0$ as $k \to \infty$, there exists $k_2 > 0$ such that inequality (1a) fails by at most $\d_k/2$ for all $k \geq k_2$ in a sufficiently small neighborhood of $\a = \frac{1}{2}$.
So we choose $k = k_2$, and then choose $\e > 0$ small enough so that $(k_2 + 1) \ee$ is less than $\frac{1}{4}$ of the normal injectivity radius of $S$.  
This completes the proof of Theorem \ref{thm:metric of negative curvature}.

  \hfill \framebox(5,5){}

\vskip 20pt

\section{Constructing an $A$-regular metric with negative sectional curvature}\label{Section:A-regular metric}
Following \cite{Belegradek complex}, let us first recall the definition of an A-regular metric.  
A Riemannian metric is called {\it A-regular} if there exists a sequence of positive numbers $\{ A_k \}_{k = 0}^\infty$ such that, for each $k$, the $k^{th}$ covariant derivative of the curvature tensor satisfies the relation $\| \nabla^k R \|_{C^0} < A_k$.  
Note that, as a consequence when $k=0$, we have that the sectional curvature is bounded from both above and below.
Two facts that are very relevant to our situation are that any metric on a compact manifold is A-regular, and that a locally symmetric metric is A-regular.  
So in order to modify the metric from Section \ref{Section:construction of metric} to make it A-regular, we only need to worry about a small neighborhood of the cusp(s) that were created when we drilled out $S$.  
In particular, we modify the metric of Section \ref{Section:construction of metric} over a region which is to the left of $\ae$, so that $h_\th = h_r$ and we can use the formulas from Section \ref{Section:h_th=h_r}.  

For completeness, we provide the details below for how we alter the metric from Section \ref{Section:construction of metric}.  
For this we are basically copying the beginning of Section 11 of \cite{Belegradek complex}.  
But in order to prove that this metric is A-regular with negative sectional curvature, we simply explain why Belegradek's argument in \cite{Belegradek complex} goes through in our setting virtually verbatim.  

Let $\ae$ be as in Section \ref{Section:construction of metric}, and let $\tau_e := \e e^{\ae}$.  
Notice that $0 < \tau_\e < 2\e$ since $e^{\ae}$ is clearly less than $2$.  
Then let $o_\e := \ln (\tau_\e)$ and $p_e := 2 \ln (\tau_\e)$.  
Both $o_\e$ and $p_\e$ are negative and approach $- \infty$ as $\e \to 0$.  
Also, $p_\e < o_\e = \ae + \ln(\e) << \ae$.  
Define 
	\begin{equation*}
	F(r) := \frac{1}{2} \cdot \frac{e^{\frac{r}{2}}}{\tau_\e + e^{\frac{r}{2}}} = \left[ \ln \left( \tau_\e + e^{\frac{r}{2}} \right) \right]'.
	\end{equation*}
Note that $F' > 0$, $F \in (0,\frac{1}{2})$, and $F(p_\e) = \frac{1}{4}$.  

The following is just a restatement of Proposition 11.1 of \cite{Belegradek complex}.  

\vskip 15pt

\begin{proposition}
For all $\e>0$ there exists a $C^1$ function $g$ such that
	\begin{itemize}
	\item  the function $g$ is positive and increasing.  
	\item  if $r \geq o_\e$ then $g$ coincides with the function $h$ from Proposition \ref{prop:warping h}.  In particular, $g(r) = e^{\frac{r}{2}}$ for $r \in [o_\e, o_\e + 1]$.  
	\item  we have that $g(r) = \tau_\e + e^{\frac{r}{2}}$ for $r \in (-\infty, p_\e]$.   
	\item  if $r \in [p_\e, o_\e]$ then $g$ is smooth, $\frac{g'}{g}$ is increasing, $\frac{g'}{g} \in \left[ \frac{1}{4}, \frac{1}{2} \right]$, and $\frac{g''}{g} > \left( \frac{g'}{g} \right)^2 \geq \frac{1}{16}$.  
	\end{itemize}
\end{proposition}

\vskip 10pt

As in Section \ref{Section:construction of metric} we can use Lemma \ref{lemma:warping lemma} to smooth out $g$ in arbitrarily small regions near $p_\e$ and $o_\e$ while controlling $g$, $g'$, and $g''$.  
Now, let $\g$ denote the Riemannian metric constructed in Section \ref{Section:construction of metric}, but with the warping function $h$ replaced by $g$.  
We then have the following Theorem (which is Theorem 11.3 in \cite{Belegradek complex}):

\vskip 15pt

\begin{theorem}\label{thm:A-regular}
For $\e > 0$ sufficiently small the metric $\g$ is A-regular, and can be smoothed near $p_\e$ and $o_\e$ via Lemma \ref{lemma:warping lemma} so that its sectional curvatures are strictly less than $0$.
\end{theorem}

\vskip 10pt

The proof in \cite{Belegradek complex} works in the present setting as well.  
All we do now is explain why this is so in the few spots where our formulas differ.

First, Belegradek uses an induction argument to show that $\g$ is A-regular.  
His argument for the base case goes through verbatim in our setting.  
In the induction step, the covariant derivatives that arise in \cite{Belegradek complex} are slightly different than what come up here (Theorem \ref{thm:connection}).  
But we still have that $\nabla_{Y_4}Y_i = 0$ for $i = 1, 2, 3$ (and recall that for us $Y_4 = \ddr$, but in \cite{Belegradek complex} $\ddr = Y_0$), and so his arguments as to why $\nabla^{k+1} R$ is bounded still apply.  

Belegradek then shows that the sectional curvature of $\g$ is strictly less than zero.  
He shows this by breaking the region $(-\infty, o_\e + \s)$, where $\s$ is a small constant from Lemma \ref{lemma:warping lemma}, up into three different parts.  
But by applying Lemma \ref{lemma:continuity lemma} we only need to consider two regions:  $(-\infty, p_\e)$ and $(p_\e, o_\e)$.  
These regions correspond to Steps 2 and 3 in Section 11 of \cite{Belegradek complex}.  

The only relevant formulas for our sectional curvature tensor that differ from those in \cite{Belegradek complex} are for $R_{2323}$ and $R_{1423}$.  
But the bounds used in \cite{Belegradek complex} still work in our setting.  
In particular, we also have that 
	\begin{equation*}
	R_{2323} < - \left( \frac{g'}{g} \right)^2 \qquad \text{and} \qquad |R_{1423}| < \frac{v}{2h^2} \left( \ln \frac{v}{h} \right)'.
	\end{equation*}  
Thus, the arguments in \cite{Belegradek complex} show that $\g$ has negative curvature, proving Theorem \ref{thm:A-regular}.

 \hfill \framebox(5,5){}

\subsection*{Acknowledgements}  
The author is indebted to both I. Belegradek and J.F. Lafont for many helpful discussions about the work contained in this paper.
The author would also like to thank I. Chatterji and K. Wang for helpful comments pertaining to Remark \ref{rmk:corollaries}.
Lastly, the author completed some of this research while on travel that was partially funded by an AMS-Simons travel grant.

\end{document}